\documentclass[a4paper,12pt]{amsart}
\input xypic
\usepackage{amssymb}
\xyoption{all}
\newtheorem{thm}{Theorem}[section]
\newtheorem{lemma}[thm]{Lemma}
\newtheorem{cor}[thm]{Corollary}
\newtheorem{prop}[thm]{Proposition}

\theoremstyle{definition}
\newtheorem{rem}[thm]{Remark}
\newtheorem{example}[thm]{Example}

\def\hbar{\bar{h}}

\def\mapsfrom{\leftarrow\!\shortmid}
\def\Rarr#1{\buildrel #1\over \longrightarrow}

\def\iso{\buildrel \sim\over\to}

\def\GS{{\mathfrak{S}}}

\def\Gb{{\mathfrak{b}}}
\def\Gg{{\mathfrak{g}}}
\def\Gh{{\mathfrak{h}}}
\def\Gm{{\mathfrak{m}}}
\def\Gn{{\mathfrak{n}}}

\def\Ggl{{\mathfrak{gl}}}
\def\Gsl{{\mathfrak{sl}}}

\def\CA{{\mathcal{A}}}
\def\CB{{\mathcal{B}}}
\def\CC{{\mathcal{C}}}
\def\CD{{\mathcal{D}}}

\def\CF{{\mathcal{F}}}

\def\CH{{\mathcal{H}}}

\def\CL{{\mathcal{L}}}
\def\CM{{\mathcal{M}}}
\def\CN{{\mathcal{N}}}
\def\CO{{\mathcal{O}}}

\def\CX{{\mathcal{X}}}
\def\CY{{\mathcal{Y}}}

\def\BC{{\mathbf{C}}}

\def\BF{{\mathbf{F}}}
\def\BG{{\mathbf{G}}}

\def\BN{{\mathbf{N}}}

\def\BP{{\mathbf{P}}}
\def\BQ{{\mathbf{Q}}}

\def\BZ{{\mathbf{Z}}}

\def\eps{\varepsilon}

\def\br{\operatorname{br}\nolimits}
\def\Br{\operatorname{Br}\nolimits}

\def\can{{\mathrm{can}}}

\def\ch{\operatorname{ch}\nolimits}

\def\Comp{\operatorname{Comp}\nolimits}

\def\End{\operatorname{End}\nolimits}

\def\Ext{\operatorname{Ext}\nolimits}

\def\GL{\operatorname{GL}\nolimits}

\def\hd{\operatorname{hd}\nolimits}

\def\Hom{\operatorname{Hom}\nolimits}

\def\Id{\operatorname{Id}\nolimits}
\def\id{\operatorname{id}\nolimits}
\def\im{\operatorname{im}\nolimits}

\def\Ind{\operatorname{Ind}\nolimits}

\def\mMod{\operatorname{\!-mod}\nolimits}

\def\opp{{\operatorname{opp}\nolimits}}

\def\mproj{\operatorname{\!-proj}\nolimits}

\def\pr{\operatorname{pr}\nolimits}

\def\rank{\operatorname{rank}\nolimits}

\def\Res{\operatorname{Res}\nolimits}

\def\sgn{\operatorname{sgn}\nolimits}

\def\soc{\operatorname{soc}\nolimits}

\def\Tr{\operatorname{Tr}\nolimits}

\def\ie{{\em i.e.}}

\def\tX{{\tilde{X}}}

\title{Derived equivalences for symmetric groups and $\Gsl_2$-categorification}
\author{Joseph Chuang and Rapha\"el Rouquier}
\address{Joseph Chuang~:
    Department of Mathematics, University of Bristol,
    University Walk, Bristol BS8 1TW, UK.
    {\it E-mail~:~}{\tt Joseph.Chuang@bris.ac.uk}}
\address{Rapha\"el Rouquier~:
        UFR de Math\'ematiques et Institut de Math\'ematiques de Jussieu
        (CNRS UMR 7586), Universit\'e Paris 7, 2 place Jussieu,
        75251 Paris Cedex 05, FRANCE.\newline
        {\it E-mail~:~}{\tt rouquier@math.jussieu.fr}}
\date\today

\begin{document}

\begin{abstract}
We define and study $\Gsl_2$-categorifications on abelian categories.
We show in particular that there is a self-derived (even homotopy)
equivalence categorifying the adjoint action of the simple reflection.
We construct categorifications for blocks of symmetric groups
and deduce that two blocks are splendidly Rickard equivalent whenever they have
isomorphic defect groups
and we show that this implies Brou\'e's abelian
defect group conjecture for symmetric groups. We give similar results
for general linear groups over finite fields.
The constructions extend
to cyclotomic Hecke algebras. We also construct categorifications
for category $\CO$ of $\Ggl_n(\BC)$ and for
rational representations of general linear groups over $\bar{\BF}_p$, where
we deduce that two blocks corresponding
to weights with the same stabilizer under
the dot action of the affine Weyl group
have equivalent derived (and homotopy)
categories, as conjectured by Rickard.
\end{abstract}

\maketitle
\tableofcontents

\section{Introduction}
The aim of this paper is to show that two blocks of symmetric groups
with isomorphic defect groups have equivalent derived categories. We deduce
in particular that Brou\'e's abelian defect group conjecture holds for
symmetric groups. We prove similar results for general linear groups
over finite fields and cyclotomic Hecke algebras.

\smallskip
Recall that there is an action of $\hat{\Gsl}_p$ on the sum of
Grothendieck groups of categories of $k\GS_n$-modules, for $n\ge 0$,
where $k$ is a field of characteristic $p$. The action of the generators
$e_i$ and $f_i$ come from exact functors between modules
(``$i$-induction'' and ``$i$-restriction''). The adjoint action
of the simple reflections of the affine Weyl group can be
categorified as functors between derived categories, following Rickard.
The key point is to show that these functors are invertible,
since two blocks have isomorphic defect groups if and only if they are in
the same affine Weyl group orbit.
This involves only an $\Gsl_2$-action and we solve the problem in a more
general framework.

\smallskip
We develop a notion of $\Gsl_2$-{\em categorification} on an
abelian category. This involves the data of adjoint exact functors
$E$ and $F$ inducing an $\Gsl_2$-action on the Grothendieck group and
the data of endomorphisms $X$ of $E$ and $T$ of $E^2$ satisfying
the defining relations of (degenerate) affine Hecke algebras.

Our main Theorem is a proof that the categorification
$\Theta$ of the simple reflection is a self-equivalence at the level
of derived (and homotopy) categories.
We achieve this in two steps. First, we show
that there is a minimal categorification of string (=simple) modules coming
from certain quotients of (degenerate) affine Hecke algebras~: this reduces
the proof of invertibility of $\Theta$ to the case of the
minimal categorification. There, $\Theta$
becomes (up to shift) a self-equivalence of the abelian category.

\medskip
Let us now describe in more detail the structure of this article.
The first part \S\ref{secaffine} is devoted to the study of
(degenerate) affine Hecke algebras of type $A$ completed at a maximal ideal
corresponding to a totally ramified central character.
We construct (in \S\ref{totram})
explicit decompositions of tensor products of
ideals which we later translate into isomorphisms of functors.
In \S\ref{quo}, we introduce certain quotients, that turn out to
be Morita equivalent to cohomology rings of Grassmannians.
Part \S\ref{reminders} recalls elementary results on
adjunctions and on representations of $\Gsl_2$.

Part \S\ref{seccat} is
devoted to the definition and study of $\Gsl_2$-categorifications. We first
define a weak version (\S\ref{weak}),
with functors $E$ and $F$ satisfying $\Gsl_2$-relations
in the Grothendieck group. This is enough to get filtrations of the
category and to introduce a class of objects that control the abelian
category. Then, in \S\ref{subseccat},
we introduce the extra data of $X$ and $T$ which give
the genuine $\Gsl_2$-categorifications. This provides actions
of (degenerate) affine Hecke algebras on powers of $E$ and $F$.
This leads immediately to two constructions of divided powers of $E$ and $F$.
In order to study $\Gsl_2$-categorifications, we introduce in \S\ref{mincat}
``minimal''
categorifications of the simple $\Gsl_2$-representations, based on
the quotients introduced in \S\ref{quo}. A key construction
(\S \ref{reducetominimal})
is a functor from such a minimal categorification to a given categorification,
that allows to reduce part of the study of an arbitrary
$\Gsl_2$-categorification to this minimal case, where explicit computations
can be carried out. This corresponds to the decomposition of the
$\Gsl_2$-representation on $K_0$ into a direct sum of irreducible
representations.
We use this in \S\ref{decompo}
to prove a categorified version of the relation $[e,f]=h$ and deduce
a construction of categorifications on the module category of the endomorphism
ring of ``stable'' objects in a given categorification.

Part \S\ref{secrefl} is devoted to the categorification of the simple
reflection of the Weyl group. In \S\ref{Rickardcplx}, we construct a
complex of functors categorifying this reflection, following Rickard.
The main result is Theorem \ref{mainthm} in part 
\S\ref{dersimple}, which shows that this complex induces a self-equivalence of
the homotopy and of the derived category. The key step in the proof for the
derived category is
the case of a minimal categorification, where we show that the complex
has homology concentrated in one degree (\S\ref{equivmin}). The case
of the homotopy category is reduced to the derived category thanks to
the constructions of \S\ref{decompo}.

In part \S\ref{secappl}, we study various examples. We define (in
\S\ref{secsymm})
$\Gsl_2$-categorifications on representations of symmetric groups and
deduce derived and even splendid Rickard equivalences.
We deduce a proof of Brou\'e's abelian defect group
conjecture for blocks of symmetric groups. We give similar
constructions for cyclotomic Hecke algebras (\S\ref{cyclo}) and for general
linear groups over a finite field in the non-defining characteristic case
(\S\ref{GLq}) for which
we also deduce the validity of Brou\'e's abelian defect group conjecture.
We also construct $\Gsl_2$-categorifications on category $\CO$ for
$\Ggl_n$ (\S\ref{CatO}) and on rational representations of $\GL_n$ over an
algebraically
closed field of characteristic $p>0$ (\S\ref{ratrep}).
This answers in particular the $\GL$ case of a conjecture
of Rickard on blocks corresponding to weights with same stabilizers under
the dot action of the affine Weyl group. We also explain similar
constructions for $q$-Schur algebras (\S\ref{qSchur}) and provide morphisms of
categorifications relating the previous constructions. A special
role is played by the endomorphism $X$, which takes various incarnations~:
the Casimir in the rational representation case and the Jucys-Murphy elements
in the Hecke algebra case. In the case of the general linear groups
over a finite field, our construction seems to be new.
Our last section (\S\ref{realmin}) provides various realizations of minimal
categorifications, including one coming from the geometry of
Grassmannian varieties.

\smallskip
Our general approach is inspired by
\cite{LLT}, \cite{Ar1}, \cite{Gr}, and \cite{GrVa}
(cf \cite[\S 3.3]{Rou3}), and our strategy for
proving the invertibility of $\Theta$ is reminiscent of \cite{DeLu,CaRi}.

\smallskip
In a work in progress, we study the braid relations between the
categorifications of the simple reflections, in the more general
framework of categorifications of Kac-Moody algebras and in relation
with Nakajima's quiver variety constructions.

\smallskip
The first author was supported in this research 
by the Nuffield Foundation (NAL/00352/G)
and the EPSRC (GR/R91151/01).

\section{Notations}
Given an algebra $A$, we denote by $A^\opp$ the opposite algebra. We
denote by $A\mMod$ the category of finitely generated $A$-modules. Given
an abelian category $\CA$, we denote
by $\CA\mproj$ the category of projective objects of $\CA$.

Let $\CC$ be an additive category. We denote by $\Comp(\CC)$ the category
of complexes of objects of $\CC$ and by $K(\CC)$ the corresponding
homotopy category.

Given an object $M$ in an abelian category, we denote by $\soc(M)$ 
(resp. $\hd(M)$) the socle (resp. the head) of $M$, \ie, the largest 
semi-simple subobject (resp. quotient) of $M$, when this exists.

We denote by $K_0(\CA)$ the Grothendieck group of an exact category $\CA$.

Given a functor $F$, we write sometimes $F$ for the identity endomorphism
$\mathbf{1}_F$ of $F$.

\section{Affine Hecke algebras}
\label{secaffine}
\subsection{Definitions}
Let $k$ be a field and $q\in k^\times$. We define a $k$-algebra
$H_n=H_n(q)$.

\subsubsection{Non-degenerate case}
Assume $q\not=1$.
The affine Hecke algebra $H_n(q)$ is the $k$-algebra with
generators
$$T_1,\ldots,T_{n-1},X_1^{\pm 1},\ldots,X_n^{\pm 1}$$
subject to the relations
\begin{align*}
(T_i+1)(T_i-q)&=0 \\
T_iT_j  &  =T_jT_i\quad (\textrm{when } \left\vert i-j\right\vert >1) \\
T_iT_{i+1}T_i  &  =T_{i+1}T_iT_{i+1}\\
X_iX_j  &  =X_jX_i\\
X_iX_i^{-1}&=X_i^{-1}X_i=1\\
X_iT_j  &  =T_jX_i\quad (\textrm{when }  i-j\not=0,1)\\
T_iX_iT_i &=qX_{i+1}.
\end{align*}

We denote by $H_n^f(q)$ the subalgebra of $H_n(q)$ generated by 
$T_1,\ldots,T_{n-1}$. It is the Hecke algebra of the symmetric group $\GS_n$.

Let $P_n=k[X_1^{\pm 1},\ldots,X_n^{\pm 1}]$, a subalgebra of $H_n(q)$
of Laurent polynomials. We put also $P_{[i]}=k[X_i^{\pm 1}]$.

\subsubsection{Degenerate case}
Assume $q=1$.
The degenerate affine Hecke algebra $H_n(1)$ is the $k$-algebra with
generators
\[
T_{1},\ldots,T_{n-1},X_{1},\ldots,X_{n}%
\]
subject to the relations
\begin{align*}
T_{i}^{2}  &  =1\\
T_{i}T_{j}  &  =T_{j}T_{i}\quad(\textrm{when }  \left\vert i-j\right\vert >1) \\
T_{i}T_{i+1}T_{i}  &  =T_{i+1}T_{i}T_{i+1}\\
X_{i}X_{j}  &  =X_{j}X_{i}\quad\\
X_{i}T_{j}  &  =T_{j}X_{i}\quad(\textrm{when }  i-j\not= 0,1)\\
X_{i+1}T_{i}&=T_{i}X_{i}+1.
\end{align*}

Note that the degenerate affine Hecke algebra is not the specialization of 
the affine Hecke algebra.

\smallskip
We put $P_n=k[X_1,\ldots,X_n]$, a polynomial subalgebra of $H_n(1)$.
We put also $P_{[i]}=k[X_i]$.
The subalgebra $H_n^f(1)$ of $H_n(1)$ generated by $T_1,\ldots,T_{n-1}$
is the group algebra $k\GS_n$ of the symmetric group.

\subsubsection{}
We put $H_n=H_n(q)$ and
$H_n^f=H_n^f(q)$.
There is an isomorphism
$H_n\iso H_n^\opp, T_i\mapsto T_i,\ X_i\mapsto X_i$.
It allows us to switch between right and left $H_n$-modules.
There is an automorphism of $H_n$ defined by
$T_i\mapsto T_{n-i},\ X_i\mapsto \tX_{n-i+1}$,
where $\tX_i=X_i^{-1}$ if $q\not=1$ and $\tX_i=-X_i$ if $q=1$.

We denote by $l:\GS_n\to\BN$ the length function.
We put $s_i=(i,i+1)\in\GS_n$.
Given $w=s_{i_1}\cdots s_{i_r}$ a reduced decomposition of
an element $w\in\GS_n$ (\ie, $r=l(w)$), we put
$T_w=T_{s_{i_1}}\cdots T_{s_{i_r}}$.

We have $H_{n}=H_n^f\otimes P_n=P_n\otimes H_n^f$.

We have an action of $\GS_n$ on $P_n$ by permutation of the variables.
Given $p\in P_n$, we have \cite[Proposition 3.6]{Lu}
\begin{equation}
\label{Lusztig}
T_ip-s_i(p)T_i=
\begin{cases}
(q-1)(1-X_iX_{i+1}^{-1})^{-1}(p-s_i(p)) & \textrm{ if }q\not=1 \\
(X_{i+1}-X_i)^{-1}(p-s_i(p)) & \textrm{ if }q=1
\end{cases}
\end{equation}
Note that $(P_n)^{\GS_n}\subset Z(H_n)$ (this is actually an equality, 
a result of Bernstein).

\subsubsection{}
\label{cn}
Let $1$ (resp. $\sgn$) be the one-dimensional representation of $H_n^f$ given
by $T_{s_i}\mapsto q$ (resp. $T_{s_i}\mapsto -1$).
Let $\tau\in \{1,\sgn\}$.
We put
$$c_n^\tau=\sum_{w\in\GS_n} q^{-l(w)}\tau(T_w)T_w.$$
 We have $c_n^{\tau}\in Z(H_n^f)$. We have
$c^1_n=\sum_{w\in\GS_n}T_w$ and
$c^{\sgn}_n=\sum_{w\in\GS_n}(-q)^{-l(w)}T_w$, and
$c^1_n c^{\sgn}_n=c^{\sgn}_n c^1_n=0$ for $n\ge 2$.

More generally, given $1\le i\le j\le n$, we denote by
$\GS_{[i,j]}$ the symmetric group on $[i,j]=\{i,i+1,\ldots,j\}$,
we define similarly $H_{[i,j]}^f$, $H_{[i,j]}$ and
we put
$c^\tau_{[i,j]}=\sum_{w\in\GS_{[i,j]}} q^{-l(w)}\tau(T_w)T_w$.

Given $I$ a subset of $\GS_n$ we put
$c_I^\tau=\sum_{w\in I} q^{-l(w)}\tau(T_w)T_w$.
We have
$$c_n^\tau=
c_{[\GS_n/\GS_i]}^\tau c_i^\tau=
c_i^\tau c^\tau_{[\GS_i\setminus \GS_n]}$$
where $[\GS_n/\GS_i]$ (resp. $[\GS_i\setminus \GS_n]$) is the set of minimal
length representatives of right (resp. left) cosets.

Given $M$ a projective $H_n^f$-module, then 
$c^\tau_n M=\{m\in M\mid hm=\tau(h)m \textrm{ for all }h\in H_n^f\}$ and
the multiplication map $c_n^\tau H_n^f \otimes_{H_n^f}M\iso c_n^\tau M$
is an isomorphism. Given $N$ an $H_n$-module, then the canonical map
$c_n^\tau H_n^f\otimes_{H_n^f}N\iso c_n^\tau H_n\otimes_{H_n}N$ is an
isomorphism.

\subsection{Totally ramified central character}
\label{totram}
We gather here a number of properties of (degenerate) affine Hecke algebras
after completion at a maximally ramified central character. Compared to
classical results, some extra complications arise from the possibility
of $n!$ being $0$ in $k$.

\subsubsection{}
\label{en}
We fix $a\in k$, with $a\not=0$ if $q\not=1$. 
We put $x_i=X_i-a$.
Let $\Gm_n$ be the maximal ideal of $P_n$ generated by $x_1,\ldots,x_n$ and let
$\Gn_n=(\Gm_n)^{\GS_n}$.

Let $e_m(x_1,\ldots,x_n)=\sum_{1\le i_1<\cdots<i_m\le n} x_{i_1}\cdots x_{i_m}
\in P_n^{\GS_n}$ be the $m$-th elementary symmetric function.
Then, $x_n^n=\sum_{i=0}^{n-1} (-1)^{n+i+1} x_n^ie_{n-i}(x_1,\ldots,x_n)$.
So, $x_n^l\in \bigoplus_{i=0}^{n-1} x_n^i \Gn_n$ for
$l\ge n$. Via Galois theory,
we deduce that 
$P_n^{\GS_{n-1}}=\bigoplus_{i=0}^{n-1} x_n^i P_n^{\GS_n}$.
Using that the multiplication map
$P_j^{\GS_j}\otimes P_{[j+1,n]}\iso P_n^{\GS_j}$ is an isomorphism, we deduce
by induction that
\begin{equation}
\label{decfix}
P_n^{\GS_{r}}=
\bigoplus_{0\le a_i< r+i} x_{r+1}^{a_1}\cdots x_n^{a_{n-r}}P_n^{\GS_n}.
\end{equation}

\subsubsection{}
\label{kato}
We denote by $\widehat{P_n^{\GS_n}}$ the completion of $P_n^{\GS_n}$ at $\Gn_n$,
and we put $\hat{P}_n=P_n\otimes_{P_n^{\GS_n}}\widehat{P_n^{\GS_n}}$ and
$\hat{H}_n=H_n\otimes_{P_n^{\GS_n}}\widehat{P_n^{\GS_n}}$. The canonical
map $\hat{P}_n^{\GS_n}\iso \widehat{P_n^{\GS_n}}$ is an isomorphism,
since $\widehat{P_n^{\GS_n}}$ is flat over $P_n^{\GS_n}$.

We denote by $\CN_n$ the category of locally nilpotent
$\hat{H}_n$-modules, \ie,
the category of $H_n$-modules on which
$\Gn_n$ acts locally nilpotently~: an $H_n$-module $M$ is
in $\CN_n$ if for every $m\in M$, there is $i>0$ such that $\Gn_n^i m=0$.

We put $\bar{H}_n=H_n/(H_n\Gn_n)$ and
$\bar{P}_n=P_n/(P_n \Gn_n)$. The multiplication gives an
isomorphism $\bar{P}_n\otimes H_n^f\iso \bar{H}_n$.
The canonical map
$$\bigoplus_{0\le a_i<i}k x_1^{a_1}\cdots x_n^{a_n}\iso
\bar{P}_n$$
 is an isomorphism, hence $\dim_k \bar{H}_n=(n!)^2$.

The unique simple object of $\CN_n$ is
 \cite[Theorem 2.2]{Ka}
$$K_n=H_n\otimes_{P_n}P_n/\Gm_n\simeq \bar{H}_n c_n^\tau.$$
It has dimension $n!$ over $k$.
It follows that the canonical surjective map
$\bar{H}_n\to \End_k(K_n)$ is an isomorphism, hence
$\bar{H}_n$ is a simple split $k$-algebra.

Since $K_n$ is a free module over $H_n^f$, it follows that
any object of $\CN_n$ is free by restriction to $H_n^f$.
From \S\ref{cn}, we deduce that for any $M\in\CN_n$, the canonical map
$c_n^\tau H_n\otimes_{H_n}M\iso c_n^\tau M$ is an isomorphism.

\begin{rem}
We have excluded the case of the
affine Weyl group algebra (the affine Hecke algebra at $q=1$). Indeed,
in that case $K_n$ is not simple (when $n\ge 2$) and $\bar{H}_n$
is not a simple algebra. When $n=2$, we have
$\bar{H}_n\simeq \left(k[x]/(x^2)\right)\rtimes \mu_2$, where the group
$\mu_2=\{\pm 1\}$ acts on $x$ by multiplication.
\end{rem}

\subsubsection{}
Let $f:M\to N$ be a morphism of finitely generated $\hat{P}_n^{\GS_n}$-modules.
Then, $f$ is surjective if and only if
$f\otimes_{\hat{P}_n^{\GS_n}}\hat{P}_n^{\GS_n}/\hat{\Gn}_n$ is surjective.

\begin{lemma}
\label{decomp0}
We have isomorphisms
$$\hat{H}_nc_n^\tau \otimes_k \bigoplus_{i=0}^{n-1}x_n^i k
\xrightarrow[\sim]{\can}
\hat{H}_nc_n^\tau\otimes_{\hat{P}_n^{\GS_n}}\hat{P}_n^{\GS_{n-1}}
\xrightarrow[\sim]{\mathrm{mult}} \hat{H}_nc^\tau_{n-1}.$$
\end{lemma}

\begin{proof}
The first isomorphism follows from the decomposition of
$\hat{P}_n^{\GS_{n-1}}$ in (\ref{decfix}).

Let us now study the second map.
Note that both terms are free $\hat{P}_n^{\GS_n}$-modules of rank
$n\cdot n!$, since
$\hat{H}_nc^\tau_{n-1}\simeq \hat{P}_n\otimes H_n^f c^\tau_{n-1}$.
Consequently, it suffices to show that the map is surjective.
Thanks to the remark above, it is enough to check surjectivity after
applying $-\otimes_{\hat{P}_n^{\GS_n}} \hat{P}_n^{\GS_n}/\hat{\Gn}_n$.

Note that the canonical surjective map
$k[x_n]\to P_n^{\GS_{n-1}}\otimes_{P_n^{\GS_n}}P_n^{\GS_n}/\Gn_n$
factors through $k[x_n]/(x_n^n)$ (cf \S\ref{en}).
So, we have to show that the multiplication map
$f:\bar{H}_nc_n^\tau\otimes k[x_n]/(x_n^n)\to \bar{H}_nc_{n-1}^\tau$
is surjective.
This is a morphism of $(\bar{H}_n,k[x_n]/(x_n^n))$-bimodules.
The elements $c_n^\tau,c_n^\tau x_n,\ldots,c_n^\tau x_n^{n-1}$ of
$\bar{H}_n$ are linearly independent, hence
the image of $f$ is a faithful $(k[x_n]/(x_n^n))$-module. It follows
that $f$ is injective, since $\bar{H}_nc_n^\tau$ is a simple $\bar{H}_n$-module.
Now, $\dim_k \bar{H}_nc_{n-1}^\tau=n\cdot n!$,
hence $f$ is an isomorphism.
\end{proof}

\smallskip
Let $M$ be a $k\GS_n$-module. We put
$\Lambda^{\GS_n}M=M/(\sum_{0<i<n}M^{s_i})$. If $n!\in k^\times$,
then $\Lambda^{\GS_n}M$ is the largest quotient of $M$ on which $\GS_n$
acts via the sign character. Note that given a vector space $V$, then
$\Lambda^{\GS_n}(V^{\otimes n})=\Lambda^n V$.

\begin{prop}
\label{decomp1}
Let $\{\tau,\tau'\}=\{1,\sgn\}$ and $r\le n$.
We have isomorphisms
$$\hat{H}_nc_n^\tau \otimes_k \bigoplus_{0\le a_i< n-r+i}
x_{n-r+1}^{a_1}\cdots x_n^{a_r} k\xrightarrow[\sim]{\can}
\hat{H}_nc_n^\tau\otimes_{\hat{P}_n^{\GS_n}}\hat{P}_n^{\GS_{[1,n-r]}}
\xrightarrow[\sim]{\mathrm{mult}} \hat{H}_nc^\tau_{[1,n-r]}.$$
There is a commutative diagram
$$\xymatrix{
& \hat{H}_nc_n^\tau\otimes_{\hat{P}_n^{\GS_n}}\hat{P}_n^{\GS_{[1,n-r]}}
\ar@{->>}[dr]^{\ x\otimes y\mapsto xy c^{\tau'}_{[n-r+1,n]}}
\ar@{->>}[d]_{\can} \\
\hat{H}_nc_n^\tau \otimes_k {\displaystyle\bigoplus_{0\le a_1<\cdots<a_r< n}}
x_{n-r+1}^{a_1}\cdots x_n^{a_r} k\ar[r]^-\sim_-{\can} &
\hat{H}_nc_n^\tau\otimes_{\hat{P}_n^{\GS_n}}
\Lambda^{\GS_{[n-r+1,n]}}\hat{P}_n^{\GS_{[1,n-r]}}\ar@{.>}[r]^-\sim &
\hat{H}_nc^\tau_{[1,n-r]}c^{\tau'}_{[n-r+1,n]}
}$$
\end{prop}

\begin{proof}
The multiplication map $H_n\otimes_{H_{n-i}}H_{n-i}c_{n-i}^\tau\to
H_n c_{n-i}^\tau$ is an isomorphism (cf \S \ref{cn}).
It follows from Lemma \ref{decomp0} that multiplication is an
isomorphism
$$\hat{H}_nc_{n-r+1}^\tau\otimes\bigoplus_{i=0}^{n-r}x_{n-r+1}^ik\iso
\hat{H}_nc_{n-r}^\tau$$
and the first statement follows by descending induction on $r$.

\smallskip
The surjectivity of the diagonal map follows from the first statement
of the Proposition.

Let $p\in \hat{P}_n^{s_i}$. Then, $c_{[i,i+1]}^1 p=pc_{[i,i+1]}^1$.
It follows that
$c_{[i,i+1]}^\tau pc_{[i,i+1]}^{\tau'}=0$, hence
$c_n^\tau p c_{[n-r+1,n]}^{\tau'}=0$ whenever $i\ge n-r+1$.
This shows the factorization property (existence of the dotted arrow).

Note that $\Lambda^{\GS_{[n-r+1,n]}}\hat{P}_n^{\GS_{n-r}}$ is generated by
$\bigoplus_{0\le a_1<\cdots<a_r< n} x_{n-r+1}^{a_1}\cdots x_n^{a_r} k$
as a $\hat{P}_n^{\GS_n}$-module (cf (\ref{decfix})).
It follows that we have surjective maps
$$\hat{H}_nc_n^\tau \otimes_k \bigoplus_{0\le a_1<\cdots<a_r<n}
x_{n-r+1}^{a_1}\cdots x_n^{a_r} k\twoheadrightarrow
\hat{H}_nc_n^\tau\otimes_{\hat{P}_n^{\GS_n}}
\Lambda^{\GS_{[n-r+1,n]}}\hat{P}_n^{\GS_{n-r}}\twoheadrightarrow
\hat{H}_nc^\tau_{n-r}c^{\tau'}_{[n-r+1,n]}.$$
Now the first and last terms above are free $\hat{P}_n$-modules of
rank $n \choose r$, hence the maps are isomorphisms.
\end{proof}

\begin{lemma}
\label{invar}
Let $r\le n$.
We have
$c_r^\tau \hat{H}_n c_n^\tau=\hat{P}_n^{\GS_r}c_n^\tau$,
$c_n^\tau \hat{H}_n c_r^\tau=c_n^\tau\hat{P}_n^{\GS_r}$ 
and the multiplication maps
$c_n^\tau\hat{H}_n\otimes_{\hat{H}_n}\hat{H}_nc_r^\tau\iso
c_n^\tau \hat{H}_n c_r^\tau $
and
$c_r^\tau\hat{H}_n\otimes_{\hat{H}_n}\hat{H}_nc_n^\tau\iso
c_r^\tau \hat{H}_n c_n^\tau $
are isomorphisms.
\end{lemma}

\begin{proof}
We have an isomorphism $\hat{P}_n\iso \hat{H}_nc_n^\tau,\ p\mapsto pc_n^\tau$.
Let $h\in\hat{H}_n$. We have $c_n^\tau hc_n^\tau=pc_n^\tau$ for some
$p\in\hat{P}_n$.
Since $T_ic_n^\tau=\tau(T_i)c_n^\tau$, it follows that
$T_i p c_n^\tau=\tau(T_i)pc_n^\tau$. So,
$(T_ip-s_i(p)T_i)c_n^\tau=\tau(T_i)(p-s_i(p))c_n^\tau$, hence
$p-s_i(p)=0$, using the formula (\ref{Lusztig}). It follows that
$c_n^\tau\hat{H}_nc^\tau_n\subseteq \hat{P}_n^{\GS_n}c_n^\tau$.

By Proposition \ref{decomp1},
the multiplication map
$\hat{H}_nc^\tau_n\otimes_{\hat{P}_n^{\GS_n}}\hat{P}_n\iso \hat{H}_n$
is an isomorphism. So, the multiplication map
$c_n^\tau\hat{H}_nc^\tau_n\otimes_{\hat{P}_n^{\GS_n}}\hat{P}_n\iso
c_n^\tau\hat{H}_n$ is an isomorphism, hence the canonical map
$c_n^\tau\hat{H}_nc^\tau_n\otimes_{\hat{P}_n^{\GS_n}}\hat{P}_n\iso
\hat{P}_n^{\GS_n}c_n^\tau\otimes_{\hat{P}_n^{\GS_n}}\hat{P}_n$ is an
isomorphism. We deduce that
$c_n^\tau\hat{H}_nc^\tau_n=\hat{P}_n^{\GS_n}c_n^\tau$.

Similarly (replacing $n$ by $r$ above), we have
$c_n^\tau\hat{P}_r^{\GS_r}c_r^\tau=c_n^\tau\hat{P}_r^{\GS_r}$. Since
$P^{\GS_r}_n=P_r^{\GS_r}P_{[r+1,n]}$ (cf \S \ref{en}), we deduce that
$$c_n^\tau\hat{H}_nc_r^\tau=c_n^\tau\hat{P}_nc_r^\tau=
c_n^\tau\hat{P}_rc_r^\tau \hat{P}_{[r+1,n]}=
c_n^\tau\hat{P}_r^{\GS_r}\hat{P}_{[r+1,n]}=
c_n^\tau\hat{P}_n^{\GS_r}.$$

\smallskip
By Proposition \ref{decomp1},
$c_n^\tau\hat{H}_n\otimes_{\hat{H}_n}\hat{H}_nc_r^\tau$ is a free
$\hat{P}_n^{\GS_r}$-module of rank $1$.
So, the multiplication map
$c_n^\tau\hat{H}_n\otimes_{\hat{H}_n}\hat{H}_nc_r^\tau\to
c_n^\tau \hat{H}_n c_r^\tau $ is a surjective morphism between
free $\hat{P}_n^{\GS_r}$-modules of rank $1$, hence it is an isomorphism.

The cases where $c_r^\tau$ is on the left are similar.
\end{proof}

\begin{prop}
\label{equiv}
The functors 
$H_nc^\tau_n\otimes_{P_n^{\GS_n}}-$ and
$c^\tau_nH_n\otimes_{H_n}-$ are inverse equivalences of categories
between the category of $P_n^{\GS_n}$-modules
that are locally nilpotent for $\Gn_n$ and $\CN_n$.
\end{prop}

\begin{proof}
By Proposition \ref{decomp1},
the multiplication map
$\hat{H}_nc^\tau_n\otimes_{\hat{P}_n^{\GS_n}}\hat{P}_n\iso \hat{H}_n$
is an isomorphism.
It follows that the morphism of
$(\hat{H}_n,\hat{H}_n)$-bimodules
$$\hat{H}_nc_n^\tau\otimes_{\hat{P}_n^{\GS_n}} c_n^\tau\hat{H}_n\iso
\hat{H}_n,\ hc\otimes ch'\mapsto hch'$$ is an isomorphism.

Since $\hat{P}_n^{\GS_n}$ is
commutative, it follows from Lemma \ref{invar} that the 
$(\hat{P}_n^{\GS_n},\hat{P}_n^{\GS_n})$-bimodules $\hat{P}_n^{\GS_n}$
and $c_n^\tau\hat{H}_n\otimes_{\hat{H}_n}\hat{H}_nc_n^\tau$ are isomorphic.
\end{proof}

\subsection{Quotients}
\label{quo}
\subsubsection{}
\label{introHbar}
We denote by
$\bar{H}_{i,n}$ the image of $H_i$ in $\bar{H}_n$ for $0\le i\le n$.
Let $\bar{P}_{i,n}=P_i/(P_i\cap (P_n \Gn_n))$. 
We have an isomorphism
$H_i^f\otimes \bar{P}_{i,n}\xrightarrow[\sim]{\textrm{mult}} \bar{H}_{i,n}$.

Since $P_n^{\GS_{[i+1,n]}}=\bigoplus_{0\le a_l\le n-l}
x_1^{a_1}\cdots x_i^{a_i}P_n^{\GS_n}$ (cf (\ref{decfix})), we deduce that
$P_i=\bigoplus_{0\le a_l\le n-l} x_1^{a_1}\cdots x_i^{a_i}\oplus
(\Gn_nP_i\cap P_i)$ and $\Gn_n P_i\cap P_i=\Gn_n P_n\cap P_i$, hence
the canonical map 
\begin{equation}
\label{decPbar}
\bigoplus_{0\le a_l\le n-l} x_1^{a_1}\cdots x_i^{a_i}\iso \bar{P}_{i,n}
\end{equation}
is an isomorphism. We will identify such a monomial
$x_1^{a_1}\cdots x_i^{a_i}$ with its image in $\bar{P}_{i,n}$.
Note that $\dim_k \bar{P}_{i,n}=\frac{n!}{(n-i)!}$.

\smallskip

The kernel of the action of $P_i^{\GS_i}$ by right multiplication on
$\bar{H}_{i,n}c_i^\tau$ is $P_i^{\GS_i}\cap\Gn_n P_n$. By Proposition
\ref{equiv}, we have a Morita equivalence between $\bar{H}_{i,n}$ and
$Z_{i,n}=P_i^{\GS_i}/(P_i^{\GS_i}\cap\Gn_n P_n)$. Note that
$\bar{H}_{i,n}c_i^\tau$ is the unique indecomposable projective
$\bar{H}_{i,n}$-module and
$\dim_k \bar{H}_{i,n}=i!\dim_k \bar{H}_{i,n}c_i^\tau$. So,
$\dim_k Z_{i,n}=\frac{1}{(i!)^2}\dim_k
\bar{H}_{i,n}={n\choose i}$ and  $Z_{i,n}=Z(\bar{H}_{i,n})$.

We denote by $P(r,s)$ the set of partitions $\mu=(\mu_1\ge\cdots\ge\mu_r\ge 0)$
with $\mu_1\le s$. Given $\mu\in P(r,s)$, we denote by $m_\mu$
the corresponding monomial symmetric function
$$m_\mu(x_1,\ldots,x_r)=\sum_{\sigma} x_1^{\mu_{\sigma(1)}}\cdots
x_r^{\mu_{\sigma(r)}}$$
where $\sigma$ runs over left coset representatives of
$\GS_r$ modulo the stabilizer of $(\mu_1,\ldots,\mu_r)$.

The isomorphism (\ref{decPbar}) shows that the canonical map from
$\bigoplus_{\mu\in P(i,n-i)}k m_\mu(x_1,\ldots,x_i)$ to $\bar{P}_{i,n}$ is
injective, with image contained in $Z_{i,n}$.
Comparing dimensions, it follows that the canonical map
$$\bigoplus_{\mu\in P(i,n-i)}k m_\mu(x_1,\ldots,x_i)\iso Z_{i,n}$$
is an isomorphism.

Also, comparing dimensions, one sees that the canonical surjective maps
$$P_i\otimes_{P_i^{\GS_i}} Z_{i,n}\iso \bar{P}_{i,n} \textrm{ and }
H_i\otimes_{P_i^{\GS_i}}Z_{i,n}\iso \bar{H}_{i,n}
$$
are isomorphisms.

\subsubsection{}
\label{grassmann}

Let $G_{i,n}$ be the Grassmannian variety of $i$-dimensional subspaces of
$\BC^n$ and $G_n$ be the variety of complete flags in $\BC^n$.
The canonical morphism $p:G_n\to G_{i,n}$ induces an injective
morphism of algebras $p^*:H^*(G_{i,n})\to H^*(G_n)$ (cohomology is
taken with coefficients in $k$).
We identify $G_n$ with $GL_n/B$, where $B$ is the stabilizer of the
standard flag $(\BC(1,0,\ldots,0)\subset\cdots\subset \BC^n)$. Let
$L_j$ be the line bundle associated to the character of $B$ given
by the $j$-th diagonal coefficient.
We have an isomorphism $\bar{P}_n\iso H^*(G_n)$ sending
$x_j$ to the first Chern class of $L_j$. It multiplies degrees by $2$.
Now,
$p^*H^*(G_{i,n})$ coincides with the image of $P_i^{\GS_i}$ in $\bar{P}_n$.
So, we have obtained an isomorphism
$$Z_{i,n}\iso H^*(G_{i,n}).$$

Since $G_{i,n}$ is projective, smooth and connected of dimension
$i(n-i)$, Poincar\'e duality says that
the cup product $H^j(G_{i,n})\times H^{2i(n-i)-j}(G_{i,n})\to
H^{2i(n-i)}(G_{i,n})$
is a perfect pairing. Via the 
isomorphism $H^{2i(n-i)}(G_{i,n})\iso k$ given by
the fundamental class,
this provides $H^*(G_{i,n})$ with a structure of a symmetric algebra.

Note that the algebra $\bar{H}_{i,n}$ is isomorphic
to the ring of $i!\times i!$ matrices over $H^*(G_{i,n})$ and it is 
a symmetric algebra. Up to isomorphism, it is independent of $a$ and $q$.

\subsubsection{}
\label{relative}

Let $i\le j$.
We have
$$\bar{H}_{j,n}=
\bar{H}_{i,n}\otimes
\bigoplus_{\substack{w\in[\GS_i\setminus\GS_j]\\0\le a_l\le n-l}}
kx_{i+1}^{a_{i+1}}\cdots x_j^{a_j}\otimes kT_w$$
hence $\bar{H}_{j,n}$ is a free $\bar{H}_{i,n}$-module of
rank $\frac{(n-i)!j!}{(n-j)!i!}$.

\begin{lemma}
\label{restK}
The $H_i$-module $c_{[i+1,n]}^\tau K_n$ has a simple socle and head.
\end{lemma}

\begin{proof}
By Proposition \ref{decomp1},
multiplication gives an isomorphism
$$\bigoplus_{0\le a_l<l}x_{i+1}^{a_1}\cdots x_n^{a_{n-i}}\otimes
c_{[i+1,n]}^\tau H_{[i+1,n]}\iso  H_{[i+1,n]},$$
hence gives an isomorphism of $\bar{H}_{i,n}$-modules
$$\bigoplus_{0\le a_l<l}x_{i+1}^{a_1}\cdots x_n^{a_{n-i}}\otimes
c_{[i+1,n]}^\tau \bar{H}_n\iso  \bar{H}_n.$$
Since $\bar{H}_n$ is a free $\bar{H}_{i,n}$-module of rank
$\frac{(n-i)!n!}{i!}$, it follows that
hence $c_{[i+1,n]}^\tau \bar{H}_n$ is a free $\bar{H}_{i,n}$-module
of rank $\frac{n!}{i!}$.
We have $\bar{H}_{i,n}\simeq i!\cdot M$ as $\bar{H}_{i,n}$-modules,
where $M$ has a simple socle and head.
Since in addition $\bar{H}_n\simeq n! \cdot K_n$ as $\bar{H}_n$-modules,
we deduce that
$c_{[i+1,n]}^\tau K_n \simeq M$ has a simple socle and head.
\end{proof}

\begin{lemma}
\label{isoB}
Let $r\le l\le n$. We have isomorphisms
$$\xymatrix{
{\displaystyle\bigoplus_{0\le a_i\le n-i}}
 x_1^{a_1}\cdots x_{l-r}^{a_{l-r}}k\otimes
 {\displaystyle\bigoplus_{\mu\in P(r,n-l)}}m_\mu(x_{l-r+1},\ldots,x_l)k
 \ar[r]^-\sim_-{a\otimes b\mapsto abc_l^\tau}
 \ar[d]_\sim^{a\otimes b\mapsto ac_n^\tau\otimes b} &
 c_{[l-r+1,l]}^\tau \bar{H}_{l,n} c_l^\tau \\
\bar{H}_{l-r,n}c_{l-r}^\tau\otimes
 {\displaystyle\bigoplus_{\mu\in P(r,n-l)}}m_\mu(x_{l-r+1},\ldots,x_l)k &
c_{[l-r+1,l]}^\tau \bar{H}_{l,n}\otimes_{\bar{H}_{l,n}}\bar{H}_{l,n}c_l^\tau
\ar[u]_\sim^{\textrm{mult}}
}$$
\end{lemma}

\begin{proof}
Let $L=\bigoplus_{\mu\in P(r,n-l),0\le a_i\le n-i}
m_\mu(x_{l-r+1},\ldots,x_l)
x_1^{a_1}\cdots x_{l-r}^{a_{l-r}}k$.

We have 
$L \cap \Gn_n P_n=0$ (cf (\ref{decPbar})), hence the canonical map
$f:L\to P_l^{\GS_{[l-r+1,l]}}\otimes_{P_l^{\GS_l}}Z_{l,n}$
is injective. Since
$\dim_k Z_{l,n}={n\choose l}$ and
$P_l^{\GS_{[l-r+1,l]}}$ is a free $P_l^{\GS_l}$-module of
rank $\frac{l!}{r!}$, it
follows that $f$ is an isomorphism. Now, we have an isomorphism
(Lemma \ref{invar})
$$\hat{P}_l^{\GS_{[l-r+1,l]}}\iso c_{[l-r+1,l]}^\tau \hat{H}_lc_l^\tau,
\ a\mapsto a c_l^\tau.$$
Consequently, the horizontal map of the Lemma is an isomorphism.

As seen in \S\ref{introHbar}, the left vertical map is an isomorphism.
By Lemma \ref{invar}, the right vertical map is also an isomorphism.
\end{proof}

\section{Reminders}
\label{reminders}
\subsection{Adjunctions}
\label{secadj}

\subsubsection{}
Let $\CC$ and $\CC'$ be two categories.
Let $(G,G^\vee)$ be an adjoint pair of functors,
$G:\CC\to\CC'$ and $G^\vee:\CC'\to\CC$~: this is the data
of two morphisms $\eta:\Id_{\CC}\to G^\vee G$ (the unit) and
$\eps:G G^\vee\to \Id_{\CC'}$ (the counit), such that
$(\eps\mathbf{1}_G)\circ(\mathbf{1}_G\eta)=\mathbf{1}_G$ and
$(\mathbf{1}_{G^\vee}\eps)\circ(\eta\mathbf{1}_{G^\vee})=\mathbf{1}_{G^\vee}$.
We have then a canonical isomorphism functorial in $X\in\CC$ and $X'\in\CC'$
$$\gamma_G(X,X'): \Hom(GX,X')\iso \Hom(X,G^\vee X'),
\ \ f\mapsto G^\vee(f)\circ\eta(X),\ \ \eps(X')\circ G(f')\mapsfrom f'.$$
Note that the data of such a functorial
isomorphism provides a structure of adjoint pair.

\subsubsection{}
\label{dualend}
Let $(H,H^\vee)$ be an adjoint pair of functors, with $H:\CC\to\CC'$.
Let $\phi\in\Hom(G,H)$. Then, we define
$\phi^\vee:H^\vee\to G^\vee$ as the composition
$$\phi^\vee:H^\vee\xrightarrow{\eta_G\mathbf{1}_{H^\vee}}
G^\vee GH^\vee \xrightarrow{\mathbf{1}_{G^\vee}\phi\mathbf{1}_{H^\vee}}
G^\vee HH^\vee \xrightarrow{\mathbf{1}_{G^\vee}\eps_H} G^\vee.$$
This is the unique map
making the following diagram commutative, for any $X\in\CC$ and $X'\in\CC'$~:
$$\xymatrix{
\Hom(HX,X')\ar[rrr]^{\Hom(\phi(X),X')}\ar[d]^\sim_{\gamma_H(X,X')} &&&
 \Hom(GX,X') \ar[d]_\sim^{\gamma_G(X,X')} \\
\Hom(X,H^\vee X')\ar[rrr]_{\Hom(X,\phi^\vee(X'))} &&& \Hom(X,G^\vee X')
}$$
We have an isomorphism
$\Hom(G,H)\iso\Hom(H^\vee,G^\vee),\phi\mapsto\phi^\vee$.
We obtain in particular an isomorphism
of monoids $\End(G)\iso\End(G^\vee)^\opp$. Given
$f\in \End(G)$, then the following diagrams commute
$$\xymatrix{
& G^\vee G \ar[dr]^{\mathbf{1}_{G^\vee}f} \\
\Id_\CC\ar[ur]^{\eta}\ar[dr]_\eta  && G^\vee G \\
& G^\vee G \ar[ur]_{f^\vee\mathbf{1}_G}
}
\hspace{1cm}
\xymatrix{
& GG^\vee\ar[dr]^\eps \\
GG^\vee\ar[ur]^{f\mathbf{1}_{G^\vee}} \ar[dr]_{\mathbf{1}_G f^\vee} & &
 \Id_{\CC'} \\
& GG^\vee\ar[ur]_\eps
}$$

\subsubsection{}

Let now $(G_1,G_1^\vee)$ and $(G_2,G_2^\vee)$ be two pairs of
adjoint functors, with $G_1:\CC'\to\CC''$ and $G_2:\CC\to\CC'$.
The composite morphisms
$$\Id_\CC\xrightarrow{\eta_2}G_2^\vee G_2
 \xrightarrow{\mathbf{1}_{G_2^\vee}\eta_1 \mathbf{1}_{G_2}}
 G_2^\vee G_1^\vee G_1G_2
\ \textrm{ and }\ 
G_1G_2G_2^\vee G_1^\vee
 \xrightarrow{\mathbf{1}_{G_1}\eps_2 \mathbf{1}_{G_1^\vee}}G_1 G_1^\vee
\xrightarrow{\eps_1}\Id_\CC$$
give an adjoint pair $(G_1G_2,G_2^\vee G_1^\vee)$.

\subsubsection{}
\label{adjcplx}
Let $F=0\to F^r\xrightarrow{d^r}F^{r+1}\to\cdots\to F^s\to 0$ be a complex
of functors from $\CC$ to $\CC'$ (with $F^i$ in degree $i$).
This defines a functor $\Comp(\CC)\to\Comp(\CC')$ by taking total complexes.

Let $(F^i,F^{i\vee})$ be adjoint pairs for $r\le i\le s$. 
Let 
$$F^\vee=0\to F^{s\vee}\xrightarrow{(d^{s-1})\vee}\cdots\to
F^{r\vee}\to 0$$
 with $F^{i\vee}$ in degree $-i$. This complex of functors
defines a functor $\Comp(\CC')\to\Comp(\CC)$.

There is an adjunction $(F,F^\vee)$ between functors on categories of
complexes, uniquely determined by the property
that given $X\in\CC$ and $X'\in\CC'$, then 
$\gamma_F(X,X'):
\Hom_{\Comp(\CC')}(FX,X')\iso \Hom_{\Comp(\CC)}(X,F^\vee X')$
is the restriction of
$$\sum_i \gamma_{F^i}(X,X'):
\bigoplus_i\Hom_{\CC'}(F^iX,X')\iso \bigoplus_i\Hom_{\CC}(X,F^{i\vee}X').$$

This extends to the case where $F$ is unbounded, under the assumption that
for any $X\in\CC$, then $F^r(X)=0$ for $|r|\gg 0$ and for 
any $X'\in\CC'$, then $F^{r\vee}(X')=0$ for $|r|\gg 0$.

\subsubsection{}
Assume $\CC$ and $\CC'$ are abelian categories.

Let $c\in\End(G)$. We put $cG=\im(c)$. We assume the canonical surjection
$G\to cG$ splits (\ie, $cG=eG$ for some idempotent $e\in\End(G)$).
Then, the canonical injection $c^\vee G^\vee\to G^\vee$ splits as well
(indeed, $c^\vee G^\vee=e^\vee G^\vee$).

Let $X\in\CC$, $X'\in\CC'$ and $\phi\in\Hom(cGX,X')$. There is
$\psi\in\Hom(GX,X')$ such that $\phi=\psi_{|cGX}$. We have a commutative
diagram
$$\xymatrix{
X\ar[r]^-\eta\ar[rd]_{\eta} & G^\vee GX\ar[r]^{G^\vee c} &
 G^\vee GX\ar[r]^{G^\vee \psi} & G^\vee X' \\
& G^\vee G X \ar[ur]_{c^\vee G}\ar[r]_{G^\vee \psi} & G^\vee X'\ar[ur]_{c^\vee}
}$$
It follows that there is a (unique) map
$\gamma_{cG}(X,X'):\Hom(cGX,X')\to \Hom(X,c^\vee G^\vee X')$
making the following diagram commutative
$$\xymatrix{
\Hom(GX,X')\ar[rr]_\sim^{\gamma_G(X,X')} && \Hom(X,G^\vee X') \\
\Hom(cGX,X')\ar@{.>}[rr]^{\sim}_{\gamma_{cG}(X,X')}\ar@{^{(}->}[u] &&
 \Hom(X,c^\vee G^\vee X') \ar@{^{(}->}[u]
}$$
where the vertical maps come from the canonical projection
$G\to cG$ and injection $c^\vee G^\vee \to G^\vee$.

Similarly, one shows there is a (unique) map
$\gamma'_{cG}(X,X'):\Hom(X,c^\vee G^\vee X')\to \Hom(cGX,X')$
making the following diagram commutative
$$\xymatrix{
\Hom(GX,X') && \Hom(X,G^\vee X') \ar[ll]^\sim_{\gamma_G(X,X')^{-1}} \\
\Hom(cGX,X')\ar@{^{(}->}[u] && \Hom(X,c^\vee G^\vee X')
 \ar@{.>}[ll]_\sim^{\gamma'_{cG}(X,X')} \ar@{^{(}->}[u]
}$$

The maps $\gamma_{cG}(X,X')$ and $\gamma'_{cG}(X,X')$
are inverse to each other and
they provide $(cG,c^\vee G^\vee)$ with the structure of an adjoint pair.
If $p:G\to cG$ denotes the canonical surjection, then
$p^\vee:c^\vee G^\vee\to G^\vee$ is the canonical injection.

\subsubsection{}
\label{commutadj}
Let $\CC$, $\CC'$, $\CD$ and $\CD'$ be four categories,
$G:\CC\to\CC'$, $G^\vee:\CC'\to\CC$,
$H:\CD\to\CD'$ and $H^\vee:\CD'\to\CD$, and
$(G,G^\vee)$ and $(H,H^\vee)$ be two adjoint pairs.
Let $F:\CC\to \CD$ and $F':\CC'\to\CD'$ be two fully faithful functors
and $\phi:F'G\iso HF$ be an isomorphism.

We have isomorphisms
\begin{multline*}
\Hom(GG^\vee,\Id_{\CC'})\xrightarrow[\sim]{F'}
\Hom(F'GG^\vee,F')\xrightarrow[\sim]{\Hom(\phi^{-1}\mathbf{1}_{G^\vee},F')}
\Hom(HFG^\vee,F')\\
\xrightarrow[\sim]{\gamma_H(FG^\vee,F')}
\Hom(FG^\vee,H^\vee F')
\end{multline*}
and let $\psi:FG^\vee\to H^\vee F'$ denote the image of $\eps_G$
under this sequence of isomorphisms.

Then, $\psi$ is an isomorphism and we have a commutative diagram
$$\xymatrix{
F'GG^\vee\ar[r]^-{\mathbf{1}_F\eps_G} \ar[d]_{\phi\mathbf{1}_{G^\vee}} & F' \\
HFG^\vee \ar[r]_-{\mathbf{1}_H\psi} & HH^\vee F'\ar[u]_{\eps_H\mathbf{1}_{F'}}
}$$

\subsection{Representations of $\Gsl_2$}
\label{secsl2}
We put
$$e=\begin{pmatrix}0&1\\0&0\end{pmatrix},\ 
f=\begin{pmatrix}0&0\\1&0\end{pmatrix}\textrm{ and }
h=ef-fe=\begin{pmatrix}1&0\\0&-1\end{pmatrix}.$$

We have 
$$s=\begin{pmatrix} 0&1\\-1&0\end{pmatrix}=\exp(-f)\exp(e)\exp(-f)$$
$$s^{-1}=\begin{pmatrix} 0&-1\\1&0\end{pmatrix}=\exp(f)\exp(-e)\exp(f)$$

We put $e_+=e$ and $e_-=f$.

Let $V$ be a locally finite representation of $\Gsl_2(\BQ)$
(\ie, a direct sum of finite dimensional representations).
Given $\lambda\in\BZ$, we denote by $V_\lambda$ the weight space of
$V$ for the weight $\lambda$ (\ie, the $\lambda$-eigenspace of $h$).

For $v\in V$, let $h_\pm(v)=\max\{i|e_\pm^iv\not=0\}$ and
$d(v)=h_+(v)+h_-(v)+1$.

\begin{lemma}
Assume $V$ is a direct sum of isomorphic simple $\Gsl_2(\BQ)$-modules
 of dimension $d$.

Let $v\in V_\lambda$.
Then,
\begin{itemize}
\item
$d(v)=d=1+2h_\pm \pm \lambda$
\item
$e_\mp^{(j)}e_\pm^{(j)}v= {h_\mp + j\choose j}\cdot {h_\pm \choose j} v$
for $0\le j\le h_\pm$.
\end{itemize}
\end{lemma}

\begin{lemma}
\label{simplerrefl}
Let $\lambda\in\BZ$ and $v\in V_{-\lambda}$. Then,
$$s(v)=\sum_{r=\max(0,-\lambda)}^{h_-(v)}
 \frac{(-1)^r}{r!(\lambda+r)!}e^{\lambda+r}f^r(v)\text{ and }
s^{-1}(v)=\sum_{r=\max(0,\lambda)}^{h_+(v)}\frac{(-1)^r}{r!(-\lambda+r)!}
e^{r}f^{-\lambda+r}(v).$$
\end{lemma}

In the following Lemma, we investigate bases of weight vectors with positivity
properties.

\begin{lemma}
\label{basissl2}
Let $V$ be a locally finite $\Gsl_2(\BQ)$-module.
Let $\CB$ be a basis of $V$ consisting of weight vectors
such that $\bigoplus_{b\in\CB} \BQ_{\ge 0}b$ is stable
under the actions of $e_+$ and $e_-$.
Let $\CL_\pm=\{b\in\CB|e_\mp b=0\}$ and given
$r\ge 0$, let $V^{\le r}=\bigoplus_{d(b)\le r}\BQ b$.

Then,
\begin{enumerate}
\item
Given $r\ge 0$, then $V^{\le r}$ is a submodule of $V$ isomorphic to
a sum of modules of dimension $\le r$.
\item
Given $b\in\CB$, we have $e_\pm^{h_\pm(b)}b\in \BQ_{\ge 0}\CL_\mp$.
\item
Given $b\in\CL_\pm$, there is $\alpha_b\in \BQ_{>0}$ such that
$\alpha_b^{-1}e_\pm^{h_\pm(b)}b\in\CL_\mp$ and the map
$b\mapsto \alpha_b^{-1}e_\pm^{h_\pm(b)}b$ is a bijection $\CL_\pm\iso\CL_\mp$.
\end{enumerate}
The following assertions are equivalent:
\begin{itemize}
\item[(i)]
Given $r\ge 0$, then $V^{\le r}$ is
the sum of all the simple submodules of $V$ of dimension $\le r$.
\item[(ii)]
$\{e_\pm^ib\}_{b\in\CL_\pm,0\le i\le h_\pm(b)}$ is a basis of $V$.
\item[(iii)]
$\{e_\pm^ib\}_{b\in\CL_\pm,0\le i\le h_\pm(b)}$ generates $V$.
\end{itemize}
\end{lemma}

\begin{proof}
Let $b\in\CB$. We have $eb=\sum_{c\in\CB} u_c c$ with $u_c\ge 0$. We have
$0=e^{h_+(b)}eb=\sum_c u_c e^{h_+(b)}c$ and
$e^{h_+(b)}c\in\bigoplus_{b'\in\CB}\BQ_{\ge 0}b'$, hence $e^{h_+(b)}c=0$
for all $c\in\CB$ such that $u_c\not=0$. So,
$h_+(c)\le h_+(b)$ for all $c\in\CB$ such that $u_c\not=0$.
Hence, (1) holds.

We have $e_\pm^{h_\pm(b)}b=\sum_{c\in\CB}v_cc$ with $v_c\ge 0$.
Since $\sum_{c\in\CB}v_ce_\pm c=0$ and
$e_\pm c\in\bigoplus_{b'\in\CB}\BQ_{\ge 0}b'$, it follows that
$e_\pm c=0$ for all $c$ such that $v_c\not=0$, hence (2) holds.

Let $b\in\CL_\pm$. We have
$e_\pm^{h_\pm(b)}b=\sum_{c\in\CB}v_cc$ with $v_c\ge 0$.
We have $e_\mp^{h_\pm(b)}e_\pm^{h_\pm(b)}b=\beta b$ for some $\beta>0$.
So, $\sum_{c\in\CB}v_c e_\mp^{h_\pm(b)}c=\beta b$. It follows that given
$c\in\CB$ with $v_c\not=0$, there is $\beta_c\ge 0$ with
$e_\mp^{h_\pm(b)}c=\beta_c b$. Since 
$h_\pm(c)=h_\mp(b)$, then
$e_\pm^{h_\pm(b)}e_\mp^{h_\pm(b)}c
=\beta_c e_\pm^{h_\pm(b)}b$
is a non-zero multiple of $c$, and it follows that there is a unique
$c$ such that $v_c\not=0$. This shows (3).

\medskip
Assume (i). We prove by induction on $r$ that
$\{e_\pm^ib\}_{b\in\CL_\pm,0\le i\le h_\pm(b)<r}$ is a basis of $V^{\le r}$
(this is obvious for $r=0$). Assume it holds for $r=d$.
The image of $\{b\in\CB|d(b)=d+1\}$ in $V^{\le d+1}/V^{\le d}$ is a basis.
This module is a multiple of the simple module of dimension
$d+1$ and $\{b\in\CL_\pm|d(b)=d+1\}$ maps to a basis of the
lowest (resp. highest) weight space of $V^{\le d+1}/V^{\le d}$
if $\pm=+$ (resp. $\pm=-$). It follows that 
$\{e_\pm^ib\}_{b\in\CL_\pm,0\le i\le d=h_\pm(b)}$ maps to a basis of
$V^{\le d+1}/V^{\le d}$. By induction, it follows that
$\{e_\pm^ib\}_{b\in\CL_\pm,0\le i\le h_\pm(b)\le d}$ is a basis of
$V^{\le d+1}$. This proves (ii).

Assume (ii). Let $v$ be a weight vector with weight $\lambda$.
We have
$v=\sum_{b\in\CL_\pm,2i= \lambda\pm h_\pm(b)}u_{b,i}e_\pm^ib$ for some
$u_{b,i}\in\BQ$. Take $s$ maximal such that there is $b\in\CL_\pm$
with $h_\pm(b)=s+i$ and $u_{b,i}\not=0$. Then,
$e_\pm^sv=\sum_{b\in\CL_\pm,i=h_{\pm b}-s}u_{b,i}e_\pm^{h_\pm(b)}b$. Since
the $e_\pm^{h_\pm(b)}b$ for $b\in\CL_\pm$ are linearly independent, 
it follows that  $e_\pm^s v\not=0$, hence $s\le h_+(v)$. So,
if $d(v)<r$, then $h_\pm(b)<r$ for all $b$ such that $u_{b,i}\not=0$.
We deduce that (i) holds.

The equivalence of (ii) and (iii) is an elementary fact of representation
theory of $\Gsl_2(\BQ)$.
\end{proof}

\section{$\Gsl_2$-categorification}
\label{seccat}
\subsection{Weak categorifications}
\label{weak}
\subsubsection{}
Let $\CA$ be an artinian and noetherian $k$-linear abelian category
with the property that the endomorphism ring of any simple object is $k$
(\ie, every object
of $\CA$ is a successive extension of finitely many simple objects and
the endomorphism ring of a simple object is $k$).

A {\em weak $\Gsl_2$-categorification} is the data of
an adjoint pair $(E,F)$ of exact endo-functors of $\CA$ such that
\begin{itemize}
\item
the action of $e=[E]$ and $f=[F]$ on $V=\BQ\otimes K_0(\CA)$ gives a
locally finite
$\Gsl_2$-representation
\item
the classes of the simple objects of $\CA$ are weight vectors
\item
$F$ is isomorphic to a left adjoint of $E$.
\end{itemize}

We denote by $\eps:EF\to \Id$ and $\eta:\Id\to FE$ the (fixed)
counit and unit of the pair $(E,F)$. We don't fix an adjunction
between $F$ and $E$.

\begin{rem}
Assume $\CA=A\mMod$ for a finite dimensional $k$-algebra $A$.
The requirement that $E$ and $F$ induce an $\Gsl_2$-action
on $K_0(\CA)$ is equivalent to the same condition for
$K_0(\CA\mproj)$. Furthermore, the perfect pairing
$K_0(\CA\mproj)\times K_0(\CA)\to\BZ,\ ([P],[S])\mapsto
\dim_k\Hom_\CA(P,S)$ induces an isomorphism of $\Gsl_2$-modules
between $K_0(\CA)$ and the dual of $K_0(\CA\mproj)$.
\end{rem}

\begin{rem}
A crucial case of application will be $\CA=A\mMod$, where $A$
is a symmetric algebra. In that case, the choice of an adjunction
$(E,F)$ determines an adjunction $(F,E)$.
\end{rem}

We put $E_+=E$ and $E_-=F$.
By the weight space of an object of $\CA$, we will mean the weight space
of its class (whenever this is meaningful).

\smallskip
Note that the opposite category $\CA^\opp$ also carries a weak
$\Gsl_2$-categorification.

Fixing an isomorphism between $F$ and a left adjoint to $E$ gives another
weak categorification, obtained by swapping $E$ and $F$. We call it the
{\em dual} weak categorification.

The {\em trivial} weak $\Gsl_2$-categorification on $\CA$ is the one given by
$E=F=0$.

\subsubsection{}
Let $\CA$ and $\CA'$ be two weak $\Gsl_2$-categorifications. A {\em morphism of
weak $\Gsl_2$-categorifications} from $\CA'$ to $\CA$ is the data of a
functor $R:\CA'\to\CA$ and of isomorphisms of functors
$\zeta_\pm:RE_\pm'\iso E_\pm R$ such that the following diagram commutes

\begin{equation}
\label{commweak}
\xymatrix{
RF'\ar[rr]^{\zeta_-}\ar[d]_{\eta RF'} && FR \\
FERF'\ar[rr]_{F\zeta_+^{-1}F'} && FRE'F'\ar[u]_{FR\eps'}
}
\end{equation}

Note that $\zeta_+$ determines $\zeta_-$, and conversely (using a
commutative diagram equivalent to the one above).

\begin{lemma}
\label{morphismweak}
The commutativity of diagram (\ref{commweak}) is
equivalent to the commutativity of either of the following two diagrams
$$
\xymatrix{
&R\ar[dl]_{R\eta'}\ar[dr]^{\eta R} \\
RF'E' \ar[r]_{\zeta_-E'}^\sim & FRE'\ar[r]_{F\zeta_+}^\sim & FER
}
\hspace{1cm}
\xymatrix{
& R \\
RE'F'\ar[ur]^{R\eps'}\ar[r]_{\zeta_+F'}^\sim & ERF'\ar[r]_{E\zeta_-}^\sim & 
 EFR\ar[ul]_{\eps R}
}$$
\end{lemma}

\begin{proof}
Let us assume diagram (\ref{commweak}) is commutative.
We have a commutative diagram
$$\xymatrix{
R\ar[rr]^-{\eta R} \ar[d]_-{R\eta'} && FER \ar[rr]^-{F\zeta_+^{-1}}
 \ar[d]_-{FER\eta'} && FRE' \ar[drr]^-{\id}\ar[d]_-{FRE'\eta'} \\
 RF'E'\ar[rr]_-{\eta RF'E'}\ar@/_4pc/[rrrrrr]_{\zeta_-E'} && FERF'E'\ar[rr]_-{F\zeta_+^{-1}F'E'} && FRE'F'E'
 \ar[rr]_-{FR\eps'E'} && FRE' 
}$$
This shows the commutativity of the first diagram of the Lemma.
The proof of commutativity of the second diagram is similar.

Let us now assume the first diagram of the Lemma is commutative.
We have a commutative diagram
$$\xymatrix{
RF'\ar[r]^{\id} \ar[dr]_{R\eta'F'}\ar[dd]_{\eta RF'} &
 RF' \ar[r]^{\zeta_-} & FR \\
& RF'E'F' \ar[u]_{RF'\eps'}\ar[dr]_{\zeta_-E'F'} \\
FERF'\ar[rr]_{F\zeta_+^{-1}F'} && FRE'F'\ar[uu]_{FR\eps'}
}$$
So, diagram (\ref{commweak}) is commutative.
The case of the second diagram is similar.
\end{proof}

Note that $R$ induces
a morphism of $\Gsl_2$-modules
$K_0(\CA'\mproj)\to K_0(\CA)$.

\begin{rem}
Let $\CA'$ be a full abelian subcategory of $\CA$ stable under
subobjects, quotients, and stable under $E$ and $F$.
Then, the canonical functor $\CA'\to\CA$ is a morphism of weak
$\Gsl_2$-categorifications.
\end{rem}

\subsubsection{}
We fix now a weak $\Gsl_2$-categorification on $\CA$ and we investigate
the structure of $\CA$.

\begin{prop}
\label{blocks}
Let $V_\lambda$ be a weight space of $V$.
Let $\CA_\lambda$ be the full subcategory of $\CA$
of objects whose class is in $V_\lambda$. Then,
$\CA=\bigoplus_\lambda \CA_\lambda$.
So, the class of an indecomposable object of $\CA$ is a weight vector.
\end{prop}

\begin{proof}
Let $M$ be an object of $\CA$ with exactly two composition factors $S_1$ and $S_2$.
Assume $S_1$ and $S_2$ are in different weight spaces. Then, there is $\eps\in\{\pm\}$
and $\{i,j\}=\{1,2\}$ such that $h_\eps(S_i)>h_\eps(S_j)$. Let
$r=h_\eps(S_i)$.
We have $E_\eps^rM\iso E_\eps^r S_i\not=0$, hence
all the composition factors of $E_{-\eps}^r E_\eps^r M$ are in the same weight space as
$S_i$.
Now,
$$\Hom(E_{-\eps}^rE_\eps^r M,M)\simeq \Hom(E_\eps^r M,E_\eps^r M)\simeq
\Hom(M,E_{-\eps}^rE_\eps^r M)$$
 and these spaces are not zero. It follows that
$M$ has a non-zero simple quotient and a non-zero simple submodule in the same
weight space as $S_i$. So, $S_i$ is both a submodule and a quotient of $M$,
hence $M\simeq S_1\oplus S_2$.

We have shown that $\Ext^1(S,T)=0$ whenever $S$ and $T$ are simple objects
in different weight spaces. The proposition follows.
\end{proof}

Let $\CB$ be the set of classes of simple objects of $\CA$. This gives
a basis of $V$ and we can apply Lemma \ref{basissl2}.

We have a categorification of the fact that a locally finite $\Gsl_2$-module
is an increasing union of finite dimensional $\Gsl_2$-modules:
\begin{prop}
Let $M$ be an object of $\CA$. Then, there is a Serre subcategory $\CA'$
of $\CA$ stable under $E$ and $F$, containing $M$ and such that
$K_0(\CA')$ is finite dimensional.
\end{prop}

\begin{proof}
Let $I$ be the set of isomorphism classes of
simple objects of $\CA$ that arise as composition
factors of $E^iF^jM$ for some $i,j$. Since $K_0(\CA)$
is a locally finite $\Gsl_2$-module, then $E^iF^jM=0$ for $i,j\gg 0$,
hence $I$ is finite. Now, the Serre subcategory $\CA'$ generated by
the objects of $I$ satisfies the requirement.
\end{proof}

We have a (weak) generation result for $D^b(\CA)$~:

\begin{lemma}
\label{generation}
Let $C\in D^b(\CA)$ such that $\Hom_{D^b(\CA)}(E^iT,C[j])=0$
for all $i\ge 0$, $j\in\BZ$ and $T$ simple object of $\CA$ such that
$FT=0$. Then, $C=0$.
\end{lemma}

\begin{proof}
Assume $C\not=0$.
Take $n$ minimal such that $H^n(C)\not=0$ and $S$ simple such that
$\Hom(S,H^nC)\not=0$. Let $i=h_-(S)$ and let $T$ be a simple submodule
of $F^iS$. Then,
$$\Hom(E^iT,S)\simeq\Hom(T,F^iS)\not=0.$$
So, $\Hom_{D(\CA)}(E^iT,C[n])\not=0$ and we are done, since
$FT=0$.
\end{proof}

There is an obvious analog of Lemma \ref{generation} using
$\Hom(C[j],F^iT)$ with $ET=0$. Since $E$ is also a right adjoint of $F$,
there are similar statements with $E$ and $F$ swapped.

\begin{prop}
\label{vanishingfunctor}
Let $\CA'$ be an abelian category and $G$ be a complex of
exact functors from $\CA$ to $\CA'$ that have exact right adjoints.
We assume that for any $M\in\CA$ (resp. $N\in\CA'$), then
$G^i(M)=0$ (resp. $G^{i\vee}(N)=0$) for $|r|\gg 0$.

Assume $G(E^iT)$ is acyclic
for all $i\ge 0$ and $T$ simple object of $\CA$ such that
$FT=0$.
Then, $G(C)$ is acyclic for all $C\in\Comp^b(\CA)$.
\end{prop}

\begin{proof}
Consider the right adjoint complex $G^\vee$ to $G$ (cf \S\ref{adjcplx}).
We have an isomorphism 
$$\Hom_{D^b(\CA)}(C,G^\vee G(D))\simeq \Hom_{D^b(\CA')}(G(C),G(D))$$
for any $C,D\in D^b(\CA)$. These spaces vanish for
$C=E^iT$ as in the Proposition.
By Lemma \ref{generation}, they vanish for all $C$. The case
$C=D$ shows that $G(D)$ is $0$ in $D^b(\CA')$.
\end{proof}

\begin{rem}
Let $\CF$ be the smallest full subcategory of $\CA$ closed under extensions
and direct summands
and containing $E^iT$ for all $i\ge 0$ and $T$ simple object of $\CA$
such that $FT=0$.
Then, in general, not every projective object of
$\CA$ is in $\CF$ (cf the case of $\GS_3$ and $p=3$ in \S\ref{secsymm}).
On the other hand, if the representation $K_0(\CA)$ is isotypic, then one
shows that every object of $\CA$ is a quotient of an object of $\CF$ and
in particular the projective objects of $\CA$ are in $\CF$.
\end{rem}

Let $V^{\le d}=\sum_{b\in\CB,d(b)\le d}\BQ b$.
Let $\CA^{\le d}$ be the full Serre subcategory of $\CA$ of objects whose class
is in $V^{\le d}$.

Lemma \ref{basissl2}(1) gives the following Proposition.

\begin{prop}
The weak $\Gsl_2$-structure on $\CA$ restricts to one on $\CA^{\le d}$ and
induces one on $\CA/\CA^{\le d}$.
\end{prop}

So, we have a filtration of $\CA$ as
$0\subseteq \CA^{\le 1}\subseteq\cdots\subseteq\CA$ compatible with
the weak $\Gsl_2$-structure. It induces the filtration
$0\subseteq V^{\le 1}\subseteq\cdots\subseteq V$. Some aspects of the study of
$\CA$ can be reduced to the study of $\CA^{\le r}/\CA^{\le r-1}$. This is
particularly interesting when $V^{\le r}/V^{\le r-1}$ is a multiple
of the $r$-dimensional simple module.

\subsubsection{}
We now investigate simple objects and the effect of $E_\pm$ on them.

\begin{lemma}
\label{socle}
Let $M$ be an object of $\CA$. Assume that $d(S)\ge r$ whenever
$S$ is a simple subobject (resp. quotient) of $M$. Then,
$d(T)\ge r$ whenever $T$ is a simple subobject (resp. quotient) of $E_\pm^i M$.
\end{lemma}

\begin{proof}
It is enough to consider the case where $M$ lies in a weight space by
Proposition \ref{blocks}.
Let $T$ be a simple subobject of $E_\pm^i M$.
Since $\Hom(E_\mp^i T,M)\simeq \Hom(T,E_\pm^i M)\not=0$,
there is $S$ a simple subobject of $M$
that is a composition factor of $E_\mp^i T$. Hence,
$d(S)\le d(E_\mp^i T)\le d(T)$.
The proof for quotients is similar.
\end{proof}

Let $\CC_r$ be the full subcategory of $\CA^{\le r}$ with objects $M$
such that
whenever $S$ is a simple submodule or a simple quotient of $M$,
then $d(S)=r$.

\begin{lemma}
\label{stabler}
The subcategory $\CC_r$ is stable under $E_\pm$.
\end{lemma}

\begin{proof}
It is enough to consider the case where $M$ lies in a single weight space by
Proposition \ref{blocks}.
Let $M\in\CC_r$ lie in a single weight space. Let $T$ be
a simple submodule of $E_\pm M$. By Lemma \ref{socle},
we have $d(T)\ge r$. On the other hand, $d(T)\le d(E_\pm M)\le d(M)$.
Hence, $d(T)=r$. Similarly, one proves the required property for simple quotients.
\end{proof}

\subsection{Categorifications}
\label{subseccat}
\subsubsection{}
An {\em $\Gsl_2$-categorification} is a weak $\Gsl_2$-categorification
with the extra data of $q\in k^\times$ and $a\in k$ with $a\not=0$ if $q\not=1$
and of $X\in\End(E)$ and $T\in\End(E^2)$ such that
\begin{itemize}
\item
$(\mathbf{1}_ET)\circ(T\mathbf{1}_E)\circ(\mathbf{1}_ET)=
(T\mathbf{1}_E)\circ(\mathbf{1}_ET)\circ(T\mathbf{1}_E)$
\item
$(T+\mathbf{1}_{E^2})\circ(T-q\mathbf{1}_{E^2})=0$
\item
$T\circ(\mathbf{1}_E X)\circ T=
\begin{cases}
qX\mathbf{1}_E & \textrm{ if }q\not=1 \\
X\mathbf{1}_E-T & \textrm{ if }q=1
\end{cases}$
\item
$X-a$ is locally nilpotent.
\end{itemize}

Let $\CA$ and $\CA'$ be two $\Gsl_2$-categorifications. A {\em morphism of
$\Gsl_2$-categorifications} from $\CA'$ to $\CA$ is a morphism
of weak $\Gsl_2$-categorifications $(R,\zeta_+,\zeta_-)$ such that
$a'=a$, $q'=q$ and the following diagrams commute
\begin{equation}
\label{comm}
\xymatrix{
RE'\ar[r]^{\zeta_+}_\sim \ar[d]_{RX'} & ER\ar[d]^{XR} \\
RE' \ar[r]_{\zeta_+}^\sim & ER}
\hspace{1cm}
\xymatrix{
RE'E'\ar[r]^{\zeta_+ E'}_\sim \ar[d]_{RT'} & ERE'\ar[r]^{E\zeta_+}_\sim &
 EER\ar[d]^{TR} \\
RE'E'\ar[r]_{\zeta_+ E'}^\sim & ERE'\ar[r]_{E\zeta_+}^\sim & EER
}
\end{equation}

\subsubsection{}
\label{actiononpowers}
We define a morphism
$\gamma_n:H_n\to \End(E^n)$ by
$$T_i\mapsto \mathbf{1}_{E^{n-i-1}}T\mathbf{1}_{E^{i-1}} \text{ and }
X_i\mapsto \mathbf{1}_{E^{n-i}}X\mathbf{1}_{E^{i-1}}.$$
With our assumptions, the $H_n$-module $\End(E^n)$ (given by
left multiplication) is in $\CN_n$.

Let $\tau\in \{1,\sgn\}$. We put $E^{(\tau, n)}=E^n c^\tau_n$, the
image of $c_n^\tau:E^n\to E^n$.
Note that the canonical map $E^n\otimes_{H_n}H_n c^\tau_n\iso
E^{(\tau, n)}$ is an isomorphism (cf \S\ref{kato}).

In the context of symmetric groups, the following Lemma is due to Puig.
It is an immediate consequence of Proposition \ref{equiv}.

\begin{lemma}
\label{divided}
The canonical map
$E^{(\tau,n)}\otimes_{P_n^{\GS_n}} c^\tau_nH_n\iso E^n$
is an isomorphism. In particular, $E^n\simeq n!\cdot E^{(\tau,n)}$
and the functor $E^{(\tau,n)}$ is a direct summand of $E^n$.
\end{lemma}

We denote by $E^{(n)}$ one of the two isomorphic functors
$E^{(1,n)}$, $E^{(\sgn,n)}$.

\smallskip
Using the adjoint pair $(E,F)$, we obtain a morphism
$H_n\to \End(F^n)^\opp$ and the definitions and
results above have counterparts
for $E$ replaced by $F$ (cf \S \ref{dualend}).

We obtain a structure of $\Gsl_2$-categorification on the dual as follows.
Put $\tX=X^{-1}$ when $q\not=1$ (resp. $\tX=-X$ when $q=1$).
We choose an adjoint pair $(F,E)$.
Using this adjoint pair,
the endomorphisms $\tX$ of $E$ and $T$ of $E^2$ provide endomorphisms
of $F$ and $F^2$. We take these as the defining endomorphisms for the dual
categorification.
We define ``$a$'' for the dual categorification
as the inverse (resp. the opposite) of $a$ for the original categorification.

\begin{rem} 
\label{scalarshift}
The scalar $a$ can be shifted~: given $\alpha\in k^\times$ when $q\not=1$
(resp. $\alpha\in k$ when $q=1$), then we can define a new
categorification by replacing $X$ by $\alpha X$ (resp. by $X+\alpha
\mathbf{1}_E$). This changes $a$ into $\alpha a$ (resp. $\alpha+a$).
So, the scalar $a$ can always be adjusted to $1$ (resp. to $0$).
\end{rem}

\begin{rem}
Assume $V$ is a multiple of the simple $2$-dimensional $\Gsl_2$-module.
Then, a weak $\Gsl_2$-categorification
consists in the data of $\CA_{-1}$ and
$\CA_1$ together with inverse equivalences $E:\CA_{-1}\iso\CA_1$ and
$F:\CA_1\iso \CA_{-1}$. 
An $\Gsl_2$-categorification is the additional data of $q,a$ and
$X\in \End(E)\simeq Z(\CA_1)$ with $X-a$ nilpotent.
\end{rem}

\begin{rem}
As soon as $V$ contains a copy of a simple $\Gsl_2$-module of dimension $3$
or more, then $a$ and $q$ are determined by $X$ and $T$.
\end{rem}

\begin{example}
\label{3dim}
Take for $V$ the three dimensional irreducible representation of $\Gsl_2$.
Let $A_{-2}=A_2=k$ and $A_0=k[x]/x^2$. We put $\CA_i=A_i\mMod$.
On $\CA_{-2}$, define $E$ to be induction $\CA_{-2}\to \CA_0$. On
$\CA_0$, $E$ is restriction $\CA_0\to \CA_2$ and $F$ is restriction
$\CA_0\to\CA_{-2}$. On $\CA_2$, then $F$ is induction
$\CA_2\to\CA_0$.
$$\xymatrix{
k \ar@<0.5ex>[rr]^-{\Ind}  && k[x]/x^2 \ar@<0.5ex>[ll]^-{\Res}
 \ar@<0.5ex>[rr]^-{\Res}  && k \ar@<0.5ex>[ll]^-{\Ind}
}$$
Let $q=1$ and $a=0$. Let $X$ be the multiplication by $x$ on
$\Res:\CA_0\to \CA_2$ and multiplication by $-x$ on $\Ind:\CA_{-2}\to\CA_0$.
Let $T\in\End_k(k[x]/x^2)$ be the automorphism swapping $1$ and $x$.
This is an $\Gsl_2$-categorification of the adjoint representation
of $\Gsl_2$. The corresponding weak categorification was constructed
in \cite{HueKho}.
\end{example}

\begin{rem}
Take for $V$ the three dimensional irreducible representation of $\Gsl_2$.
Let $A_{-2}=A_2=k[x]/x^2$ and $A_0=k$. We put $\CA_i=A_i\mMod$.
On $\CA_{-2}$, then $E$ is restriction $\CA_{-2}\to \CA_0$. On
$\CA_0$, $E$ is induction $\CA_0\to \CA_2$ and $F$ is induction
$\CA_0\to\CA_{-2}$. On $\CA_2$, then $F$ is restriction
$\CA_2\to\CA_0$.
$$\xymatrix{
k[x]/x^2 \ar@<0.5ex>[rr]^-{\Res}  && k \ar@<0.5ex>[ll]^-{\Ind}
 \ar@<0.5ex>[rr]^-{\Ind}  && k[x]/x^2 \ar@<0.5ex>[ll]^-{\Res}
}$$

This is a weak $\Gsl_2$-categorification but not an $\Gsl_2$-categorification,
since $E^2:\CA_{-2}\to\CA_2$ is $(k[x]/x^2)\otimes_k -$, which is an
indecomposable functor.
\end{rem}

\begin{rem}
Let $A_{-2}=k$, $A_0=k\times k$ and $A_{-2}=k$. We define $E$ and $F$
as the restriction and induction functors in the same way as in
Example \ref{3dim}. Then, $V$ is the direct sum of a $3$-dimensional simple
representation and a $1$-dimensional representation. 
Assume there is $X\in\End(E)$ and $T\in\End(E^2)$ giving an
$\Gsl_2$-categorification. We have $\End(E^2)=\End_k(k^2)$ and $X1_E=1_EX=
a 1_{E^2}$. But the quotient of $H_2(q)$ by the relation $X_1=X_2=a$
is zero ! So, we have a contradiction (note here it is crucial to
exclude the affine Hecke algebra at $q=1$).
So, this is a weak $\Gsl_2$-categorification but not an
 $\Gsl_2$-categorification (note that we still have $E^2\simeq E\oplus E$).
\end{rem}

\subsection{Minimal categorification}
\label{mincat}
We introduce here a categorification of the (finite dimensional) simple
$\Gsl_2$-modules.

We fix $q\in k^\times$ and $a\in k$ with $a\not=0$ if $q\not=1$.
Let $n\ge 0$ and $B_i=\bar{H}_{i,n}$ for $0\le i\le n$.

We put $\CA(n)_\lambda=B_{(\lambda+n)/2}\mMod$ and
$\CA(n)=\bigoplus_i B_i\mMod$. We put
$E=\bigoplus_{i<n} \Ind_{B_i}^{B_{i+1}}$ and
$F=\bigoplus_{i>0} \Res^{B_i}_{B_{i-1}}$. Note that the functors
$\Ind_{B_i}^{B_{i+1}}=B_{i+1}\otimes_{B_i}-$ and
$\Res_{B_i}^{B_{i+1}}=B_{i+1}\otimes_{B_{i+1}}-$ are left and right adjoint.

We have $EF(B_i)\simeq B_i\otimes_{B_{i-1}}B_i\simeq i(n-i+1)B_i$ and
$FE(B_i)\simeq B_{i+1}\simeq (i+1)(n-i)B_i$ as left $B_i$-modules
(cf \S \ref{relative}). So, 
$(ef-fe)([B_i])=(2i-n)[B_i]$.
Now, $\BQ\otimes K_0(\CA(n)_\lambda)=\BQ [B_{(\lambda+n)/2}]$,
hence $ef-fe$ acts
on $K_0(\CA(n)_\lambda)$ by $\lambda$. It follows that $e$ and $f$ induce
an action of $\Gsl_2$ on $K_0(\CA(n))$, hence we have a weak
$\Gsl_2$-categorification.

\smallskip
The image of $X_{i+1}$ in $B_{i+1}$ gives an endomorphism of
$\Ind_{B_i}^{B_{i+1}}$ by right multiplication on
$B_{i+1}$. Taking the sum over all $i$, we get an
endomorphism $X$ of $E$.
Similarly, the image of $T_{i+1}$ in $B_{i+2}$ gives an endomorphism 
of $\Ind_{B_i}^{B_{i+2}}$ and taking the sum over all $i$, we get an
endomorphism $T$ of $E^2$.

This provides an $\Gsl_2$-categorification. The 
representation on $K_0(\CA(n))$ is the simple $(n+1)$-dimensional $\Gsl_2$-module.

\subsection{Link with affine Hecke algebras}

\subsubsection{}
The following Proposition generalizes and strengthens
results of Kleshchev \cite{Kl1,Kl2} in the symmetric
groups setting and of Grojnowski and Vazirani \cite{GrVa} in the context
of cyclotomic Hecke algebras (cf \S\ref{secsymm} and \S\ref{cyclo}).

\begin{prop}
\label{hdsoc}
Let $S$ be a simple object of $\CA$, let $n=h_+(S)$ and
$i\le n$.

\begin{itemize}
\item[(a)]
$E^{(n)}S$ is simple.
\item[(b)]
The socle and head of $E^{(i)}S$ are isomorphic to a simple object $T$ of $\CA$.
We have isomorphisms of $(\CA,H_i)$-bimodules:
$\soc E^i S\simeq \hd E^i S\simeq T\otimes K_i$.
\item[(c)]
The morphism $\gamma_i(S):H_i\to\End(E^iS)$ factors through
$\bar{H}_{i,n}$ and induces an isomorphism $\bar{H}_{i,n}\iso\End(E^iS)$.
$$\xymatrix{
& H_i\ar[dl]_-{\can} \ar[dr]^-{\gamma_i(S)} \\
\bar{H}_{i,n}\ar@{.>}[rr]_-\sim & & \End(E^iS)
}$$
\item[(d)]
We have $[E^{(i)}S]-{n\choose i}[T]\in V^{\le d(T)-1}$.
\end{itemize}
The corresponding statements with $E$ replaced by $F$ and $h_+(S)$ by
$h_-(S)$ hold as well.
\end{prop}

\begin{proof}
$\bullet\ $Let us assume (a) holds. We will show that (b), (c), and (d)
follow.

We have $E^n S\simeq n!\cdot T'$ for some $T'$ simple. So, we have
$E^n S\simeq T'\otimes R$ as $(\CA,H_n)$-bimodules, where $R$ is a right
$H_n$-module in $\CN_n$. Since $\dim R=\dim K_n$, it follows that
$R\simeq K_n$.

We have $E^{n-i} \soc E^{(i)}S\subset E^{n-i}E^{(i)}S\simeq
T'\otimes K_n c_i^1$. Since $T'\otimes K_n c_i^1$ has a simple socle
 (Lemma \ref{restK}),
it follows that $E^{n-i} \soc E^{(i)}S$ is an indecomposable
$(\CA,H_{n-i})$-bimodule. If $T$ is a non-zero summand of
$\soc E^{(i)}S$, then $E^{n-i}T\not=0$ (Lemma \ref{stabler}). So,
$T=\soc E^{(i)}S$ is simple. We
have $\soc E^i S\simeq T\otimes R$ for some $H_i$-module $R$ in $\CN_i$.
Since $\dim R=i!$, it follows that $R\simeq K_i$.
The proof for the head is similar.

\smallskip
The dimension of $\End(E^{(i)}S)$ is at most the multiplicity $p$ of
$T$ as a composition factor of $E^{(i)}S$. Since $E^{(n-i)}T\not=0$,
it follows that
the dimension of $\End(E^{(i)}S)$ is at most the number of
composition factors of $E^{(n-i)}E^{(i)}S$.
We have $E^{(n-i)}E^{(i)}S\simeq {n\choose i}\cdot
T'$. So,
$\dim\End(E^{(i)}S)\le {n\choose i}$ and
$\dim \End(E^iS)\le (i!)^2 {n\choose i}=\dim \bar{H}_{i,n}$.

Since $\ker\gamma_n(S)$ is a proper ideal of $H_n$, we have
$\ker\gamma_n(S)\subset \Gn_nH_n$.
We have $\ker\gamma_i(S)\subset H_i\cap \ker\gamma_n(S)\subset H_i\cap (\Gn_nH_n)$.
So, the canonical map $H_i\to \bar{H}_{i,n}$
factors through a surjective map:
$\im\gamma_i(S)\twoheadrightarrow\bar{H}_{i,n}$.
We deduce that $\gamma_i(S)$ is surjective and
$\bar{H}_{i,n}\iso \End(E^iS)$. So, (c) holds. We deduce also
that $p={n\choose i}$ and that if $L$ is a composition factor
of $E^{(i)}S$ with $E^{(n-i)}L\not=0$, then $L\simeq T$. So,
(d) holds.
Since the simple object $\hd E^{(i)}S$ is not killed by
$E^{(n-i)}$ (Lemma \ref{stabler}), we deduce that $\hd E^{(i)}S\simeq T$.
We have now shown (b).

\smallskip
$\bullet\ $Let us show that (a) (hence (b), (c), and (d)) holds when $FS=0$.
By Lemma \ref{basissl2} (3),
we have $[E^{(n)}S]=r[T]$ for some simple object $T$ and $r\ge 1$ integer.
Since $[F^{(n)}E^{(n)}S]=[S]$, we have $r=1$, so (a) holds.

\smallskip
$\bullet\ $Let us now show (a) in general.
Let $L$ be a simple quotient
of $F^{(r)}S$, where $r=h_-(S)$. Since
$\Hom(S,E^{(r)}L)\simeq \Hom(F^{(r)}S,L)\not=0$, we deduce that
$S$ is isomorphic to a submodule of $E^{(r)}L$.
Since $FL=0$, we know by (a) that
$E^{(n)}E^{(r)}L\simeq {n+r\choose r}T$ for some simple object $T$.
So, $E^{(n)}S\simeq m T$ for some positive integer $m$.
We have $\Hom(E^{(n)}S,T)\simeq \Hom(S,F^{(n)}T)$. Since $ET=0$,
we deduce that $\soc F^{(n)}T$ is simple (we use (b) in its ``$F$'' version).
So, $\dim \Hom(S,F^{(n)}T)\le 1$, hence $m=1$ and (a) holds.
\end{proof}

\begin{cor}
The $\Gsl_2(\BQ)$-module
$V^{\le d}$ is the sum of the simple submodules of $V$ of dimension $\le d$.
\end{cor}

\begin{proof}
Let $S$ be a simple object of $\CA$ with $r=h_-(S)$.
By Proposition \ref{hdsoc} (a), $T=F^{(r)}S$ is simple. 
We deduce that
$S\simeq\soc E^{(r)}T$ by adjunction. Now, Proposition \ref{hdsoc} (d)
shows that $[E^{(r)}T]-{d(S)\choose r}[S]\in V^{\le d(S)-1}$.

We deduce by induction on $r$ that 
$\{[E^r T]\}$ generates $V$, where $T$ runs over the
isomorphism classes of simple objects killed by $F$ and $0\le r\le h_+(T)$.
The Corollary follows from Lemma \ref{basissl2}, (iii)$\Longrightarrow$(i).
\end{proof}

\begin{rem}
Let $S$ be a simple object of $\CA$ and $i\le h_+(S)$. The action of
$Z_{i,n}=Z(\bar{H}_{i,n})$ on $E^iS$ restricts to an action on
$E^{(i)}S$. Since $E^iS$ is a faithful right $\bar{H}_{i,n}$-module,
it follows from Proposition \ref{equiv} that $E^{(i)}S$ is a faithful
$Z_{i,n}$-module. Now, $\dim \End_\CA(E^{(i)}S)=\frac{1}{(i!)^2}
\dim \bar{H}_{i,n}=\dim Z_{i,n}$, hence the morphism
$Z_{i,n}\to \End_\CA(E^{(i)}S)$ is an isomorphism.
\end{rem}

Let us now continue with the
following crucial Lemma whose proof uses some of the ideas of
the proof of Proposition \ref{hdsoc}.

\begin{lemma}
\label{commut2}
Let $U$ be a simple object of $\CA$ such that $FU=0$.
Let $n=h_+(U)$, $i<n$, and $B_i=\bar{H}_{i,n}$.
The composition of
$\eta(E^iU)\otimes 1:E^iU\otimes_{B_i}B_{i+1}\to
FE^{i+1}U\otimes_{B_i}B_{i+1}$ with the action map
$FE^{i+1}U\otimes_{B_i}B_{i+1} \to FE^{i+1}U$
is an isomorphism
$$E^i U\otimes_{B_i}B_{i+1}\iso FE^{i+1}U.$$
\end{lemma}

\begin{proof}
By Proposition \ref{equiv}, it is enough to prove that the map becomes an
isomorphism
after applying $-\otimes_{B_{i+1}}B_{i+1}c_{i+1}^1$.
By (\ref{decPbar}), we have
$B_{i+1}c_{i+1}^1=\bigoplus_{a=0}^{n-i-1}\bar{P}_{i,n}x_{i+1}^ac_{i+1}^1$.
Consider the composition 
$$\phi=g\circ (f\otimes 1):E^{(i)}U\otimes \bigoplus_{a=0}^{n-i-1}kx^a\to
FE^{(i+1)}U$$
where
$f:E^{(i)}U\xrightarrow{\eta(E^{(i)}U)}
FEE^{(i)}U\xrightarrow{\mathbf{1}_F c_{[\GS_i\setminus\GS_{i+1}]}^1U}
FE^{(i+1)}U$ and
$g:FE^{(i+1)}U\otimes \bigoplus_{a=0}^{n-i-1}kx^a\to
FE^{(i+1)}U$ is given by the action on $F$.
We have to prove that $\phi$ is an isomorphism.
We have $[FE^{(i+1)}U]=(n-i)[E^{(i)}U]$, hence it suffices to prove that
$\phi$ is injective. In order to do that, one may restrict $\phi$
to a map between the socles of the objects (viewed in $\CA$).
Let $\phi_a:\soc E^{(i)}U\to FE^{(i+1)}U$ be the restriction of $\phi$
to the socle of $E^{(i)}U\otimes kx^a$.
Since $\soc(E^{(i)}U)$ is simple (Proposition \ref{hdsoc}),
the problem is to prove that
the maps $\phi_a$ for $0\le a\le n-i-1$
are linearly independent. By adjunction, it is equivalent to
prove that the maps
$$\psi_a:E\soc E^{(i)}U\xrightarrow{x^a\mathbf{1}_{\soc E^{(i)}U}}
E\soc E^{(i)}U\xrightarrow{c^1_{[\GS_i\setminus\GS_{i+1}]}U} E^{(i+1)}U$$
are linearly independent.

\smallskip
We have
$\soc E^{i+1}U\simeq S\otimes K_{i+1}$
as $(\CA,H_{i+1})$-bimodules, where
$S=\soc E^{(i+1)}U$ is simple (Proposition \ref{hdsoc}).
Consider the right
$(k[x_{i+1}]\otimes H_i)$-submodule
$L'=\Hom_\CA(S,\soc (E\soc E^iU)))$ of $L=\Hom_\CA(S,\soc E^{i+1}U)$.
We have $H_{i+1}=(H_i\otimes P_{[i+1]})H_{i+1}^f$, hence
$L=L'H_{i+1}^f$ since $L$ is a simple right
$H_{i+1}$-module.
So, $L'c_{i+1}^1=Lc_{i+1}^1$, hence
$\soc (E\soc E^iU))c_{i+1}^1=\soc E^{(i+1)}U$. In particular,
the map $E\soc E^{(i)}U\xrightarrow{c^1_{[\GS_{i+1}/\GS_i]}U} E^{(i+1)}U$
is injective, since $E\soc E^{(i)}U$ has a simple socle by Proposition
\ref{hdsoc}.

So, we are left with proving that the maps
$E\soc E^{(i)}U\xrightarrow{X^a\mathbf{1}_{\soc E^{(i)}U}}
E\soc E^{(i)}U$ are linearly independent, \ie, that
the restriction of
$\gamma_1(T):H_1\to \End_{\CA}(ET)$ to $\bigoplus_{a=0}^{n-i-1}kX_1^a$
is injective, where $T=\soc E^{(i)}U$.
Let $I$ be the kernel of $\gamma_{n-i}(T):H_{n-i}\to\End_\CA(E^{n-i}T)$.
Then, as in the proof of Proposition \ref{hdsoc}, we have
$I\subset \Gn_{n-i}H_{n-i}$. So, $\ker\gamma_1\subset H_1\cap
\Gn_{n-i}H_{n-i}$, hence
the canonical map
$\bigoplus_{a=0}^{n-i-1}kX_1^a\to \End_\CA(E^{n-i}T)$ is injective
(cf (\ref{decPbar})) and we are done.
\end{proof}

\subsubsection{}
\label{reducetominimal}
We fix $U$ a simple object of $\CA$ such that $FU=0$.
Let $n=h_+(U)$. We put $B_i=\bar{H}_{i,n}$ for $0\le i\le n$.

The canonical isomorphisms of functors
$$E(E^iU\otimes_{B_i}-)\iso E^{i+1}U\otimes_{B_i}-\iso
E^{i+1}U\otimes_{B_{i+1}}B_{i+1}\otimes_{B_i}-$$
make the following diagram commutative
$$\xymatrix{
B_{i+1}\mMod\ar[rrr]^-{E^{i+1}U\otimes_{B_{i+1}}-} &&& \CA \\
B_i\mMod\ar[rrr]_-{E^iU\otimes_{B_i}-}
\ar[u]^{B_{i+1}\otimes_{B_i}-} &&& \CA\ar[u]_E
}$$

The canonical isomorphism of functors from Lemma \ref{commut2}
$$E^i U\otimes_{B_i}B_{i+1}\otimes_{B_{i+1}}-
\iso F(E^{i+1}U\otimes_{B_{i+1}}-)$$
make the following diagram commutative
$$\xymatrix{
B_{i+1}\mMod\ar[rrr]^-{E^{i+1}U\otimes_{B_{i+1}}-}
\ar[d]_{B_{i+1}\otimes_{B_{i+1}}-} &&& \CA\ar[d]^F \\
B_i\mMod\ar[rrr]_-{E^iU\otimes_{B_i}-} &&& \CA
}$$

\begin{thm}
\label{morphismminimal}
The construction above is a morphism of 
$\Gsl_2$-categorifications $R_U:\CA(n)\to \CA$.
\end{thm}

\begin{proof}
The commutativity of diagram (\ref{commweak}) follows from the very
definition of $\zeta_-$ given by Lemma \ref{commut2}.
The commutativity of the diagram (\ref{comm}) is obvious.
\end{proof}

\begin{rem}
Let $I_n$ be the set of isomorphism classes of simple objects $U$ of
$\CA$ such that $FU=0$ and $h_+(U)=n$.
We have a morphism of $\Gsl_2$-categorifications
$$\sum_{n,U\in I_n} R_U:\bigoplus_{n,U\in I_n} \CA(n)\to\CA$$
that is not an equivalence in general but that induces an isomorphism
$$\bigoplus_{n,U\in I_n} \BQ\otimes K_0(\CA(n)\mproj)\iso
\BQ\otimes K_0(\CA)$$
giving a canonical decomposition of $\BQ\otimes K_0(\CA)$ into simple
summands.
In that sense, the categorifications $\CA(n)$ are minimal.
\end{rem}

The following Proposition is clear.

\begin{prop}
\label{mincatequiv}
Assume $\BQ\otimes K_0(\CA)$ is a simple $\Gsl_2$-module of dimension $n+1$.
Let $U$ be the unique simple object of $\CA$ with $FU=0$.

Then, $R_U:\CA(n)\to\CA$ is an equivalence of categories if and only
if $U$ is projective.
\end{prop}

Note that a categorification corresponding to an isotypic representation
needs not be isomorphic to a sum of minimal categorifications (take
for example a trivial $\Gsl_2$-representation).

\subsection{Decomposition of $[E,F]$}
\label{decompo}
\subsubsection{}

Let $\sigma:EF\to FE$ be given as the composition
$$EF\xrightarrow{\eta \mathbf{1}_{EF}}FEEF\xrightarrow{\mathbf{1}_F T
\mathbf{1}_F} FEEF\xrightarrow{\mathbf{1}_{FE}\eps} FE.$$

The following gives the categorification of the relation
$[e,f]=h$.

\begin{thm}
\label{decH}
Let $\lambda\ge 0$.
Then, we have isomorphisms
$$\sigma+\sum_{j=0}^{\lambda-1} (\mathbf{1}_F X^j)\circ\eta\ :\ 
EF\Id_{\CA_{-\lambda}} \oplus \Id_{\CA_{-\lambda}}^{\bigoplus\lambda}\iso
FE\Id_{\CA_{-\lambda}}$$
and
$$\sigma+\sum_{j=0}^{\lambda-1} \eps\circ(X^j\mathbf{1}_F)\ :\ 
EF\Id_{\CA_{\lambda}}\iso
FE\Id_{\CA_{\lambda}}\oplus \Id_{\CA_{\lambda}}^{\bigoplus\lambda}.$$
\end{thm}

\begin{proof}
By Proposition \ref{vanishingfunctor},
it is enough to check that the maps
are isomorphisms after evaluating the functors at $E^iU$, where
$i\ge 0$ and $U$ is a simple object of $\CA_{-\lambda-2i}$
(resp. of $\CA_{\lambda-2i}$) such that $FU=0$.
Thanks to Lemma \ref{morphismweak} and Theorem \ref{morphismminimal},
we can do this with
$\CA$ replaced by a minimal categorification $\CA(n)$ and this
is the content of Proposition \ref{mackey} below.
\end{proof}

In the case of cyclotomic Hecke algebras, Vazirani \cite{Va}
had shown that the values of the functors on simple objects are isomorphic.

\begin{cor}
\label{homotopyK0}
The functors $E$ and $F$ induce an action of $\Gsl_2$ on the Grothendieck
group of $\CA$, viewed as an additive category.
\end{cor}

\subsubsection{}

We put $\gamma=
\begin{cases}
(q-1)a & \text{ if }q\not=1 \\
1 & \text{ if }q=1
\end{cases}$

$$\text{and }\ m_{ij}(c)=
\begin{cases}
\sum_{j\le d_1<\cdots<d_{i-j-c}\le i-1}
T_{d_1}\cdots T_{d_{i-j-c}} & \text{ if }c<i-j\\
1 & \text{ if }c=i-j\\
0 & \text{ if }c>i-j.
\end{cases}$$

\begin{lemma}
\label{modulo1}
Let $j<i$ and $c\ge 0$.
We have
$$T_j T_{j+1}\cdots T_{i-1} x_i^c=
\gamma^c m_{ij}(c) \pmod{\Gm_i H_i}.$$
In particular, $T_j T_{j+1}\cdots T_{i-1} x_i^c\in \Gm_i H_i$
if $c>i-j$.
\end{lemma}

\begin{proof}
By (\ref{Lusztig}), we have
$$T_{i-1}x_i^c-x_{i-1}^cT_{i-1}=
\begin{cases}
(q-1)(x_i+a)(x_{i-1}^{c-1}+x_{i-1}^{c-2}x_i+\cdots+x_i^{c-1})&
 \text{ if }q\not=1 \\
x_{i-1}^{c-1}+x_{i-1}^{c-2}x_i+\cdots+x_i^{c-1}& \text{ if }q=1.
\end{cases}$$
Hence
$$T_j T_{j+1}\cdots T_{i-1} x_i^c=T_j T_{j+1}\cdots T_{i-2} x_{i-1}^cT_{i-1}+
\gamma T_j T_{j+1}\cdots T_{i-2}x_{i-1}^{c-1}\pmod{\Gm_i H_i}.$$
Since
$m_{ij}(c)=m_{i-1,j}(c-1)+m_{i-1,j}(c)T_{i-1}$,
the Lemma follows by induction.
\end{proof}

\begin{lemma}
\label{modulo2}
Let $j\le i$, $c\ge 1$ and $e=\inf(c-1,i-j)$.
Then, we have
\begin{multline*}
T_j T_{j+1}\cdots T_i x_i^c-T_jT_{j+1}\cdots T_{i-1}x_{i+1}^cT_i=\\
\alpha\left(
\gamma^{e} x_{i+1}^{c-e-1} m_{ij}(e)+
\gamma^{e-1} x_{i+1}^{c-e}m_{ij}(e-1)+\cdots+
x_{i+1}^{c-1}m_{ij}(0) \right)
\pmod{\Gm_i H_{i+1}}
\end{multline*}
where $\alpha=
\begin{cases}
(1-q)(x_{i+1}+a) & \text{ if }q\not=1 \\
-1 & \text{ if }q=1.
\end{cases}$
\end{lemma}

\begin{proof}
We have
$$T_j T_{j+1}\cdots T_i x_i^c-T_jT_{j+1}\cdots T_{i-1}x_{i+1}^cT_i=
\alpha T_j\cdots T_{i-1}(x_i^{c-1}+\cdots+x_{i+1}^{c-1})$$
and the result follows from Lemma \ref{modulo1}.
\end{proof}

The following is a Mackey decomposition for the algebras $B_i=\bar{H}_{i,n}$.

\begin{prop}
\label{mackey}
Let $i\le n/2$. Then, we have an isomorphism of $(B_i,B_i)$-bimodules
\begin{align*}
B_i\otimes_{B_{i-1}}B_i\oplus B_i^{\oplus n-2i}&\iso B_{i+1}\\
(b\otimes b',b_1,\ldots,b_{n-2i})&\mapsto bT_ib'+\sum_{j=1}^{n-2i} b_j
X_{i+1}^{j-1}.
\end{align*}
Let now $i\ge n/2$. Then, we have an isomorphism of $(B_i,B_i)$-bimodules
\begin{align*}
B_i\otimes_{B_{i-1}}B_i&\iso B_{i+1} \oplus B_i^{\oplus 2i-n}\\
b\otimes b'&\mapsto (bT_ib',bb',bX_ib',\ldots,bX_i^{2i-n-1}b').
\end{align*}
\end{prop}

\begin{proof}
Let us consider the first map.
We know already that both sides are free $B_i$-modules
of the same rank (cf \S\ref{mincat}), hence it is enough to show surjectivity.

Let $M=(P_i/\Gm_i)\otimes_{P_i}B_{i+1}$.
This is a right $B_i$-module quotient of $B_{i+1}$.
Let $L$ be the right $B_i$-submodule of $M$ generated by
$B_iT_i+\sum_{j=0}^{n-2i-1} X_{i+1}^{j}k$.
The first isomorphism will follow from the proof that $M=L$.
From now on, all elements are viewed in $M$.

\smallskip
We have
$$x_{i+1}^{n-i}=\sum_{j=0}^{n-i-1} (-1)^{n-i-1+j}
x_{i+1}^j e_{n-i-j}(x_{i+1},\ldots,x_n).$$
Given $r\ge 2$ and $j\le n-i-1$, we have
$$e_{n-i-j}(x_r,\ldots,x_n)=
e_{n-i-j}(x_{r-1},x_r,\ldots,x_n)-x_{r-1}
e_{n-i-j-1}(x_r,\ldots,x_n).$$
Since $e_{n-i-j}(x_1,\ldots,x_n)=0$,
it follows that
$e_{n-i-j}(x_{i+1},\ldots,x_n)=0$.
So, we have $x_{i+1}^{n-i}=0$.

\smallskip
Take $1\le r\le i$. Then, $r\le n-i$ and
we have (Lemma \ref{modulo2})
\begin{multline*}
T_{i-r+1}T_{i-r+2}\cdots T_ix_i^{n-i}=\\
x_{i+1}^{n-i}T_{i-r+1}\cdots T_i+
\alpha\left(
\gamma^{r-1}x_{i+1}^{n-i-r}+\gamma^{r-2}x_{i+1}^{n-i-r+1}m_{i,i-r+1}(r-2)+
\cdots+x_{i+1}^{n-i-1}m_{i,i-r+1}(0)\right).
\end{multline*}
So, 
$$T_{i-r+1}T_{i-r+2}\cdots T_ix_i^{n-i}+\alpha\gamma^{r-1} x_{i+1}^{n-i-r}
\in \sum_{j\ge 0} x_{i+1}^{n-i-r+1+j}H_i.$$
Since $x_{i+1}^{n-i}=0$, we deduce by induction on $r$
that $x_{i+1}^{n-i-r}\in L$ for $1\le r\le i$. So,
$x_{i+1}^a\in L$ for all $a\ge 0$.
We deduce from Lemma \ref{modulo2} that
$x_{i+1}^a T_j\cdots T_i\in L$ for all $1\le j\le i$ and $a\ge 0$.
Since $B_{i+1}=\bigoplus_{0\le a\le n-i-1,w\in [\GS_{i+1}/\GS_{i}]}
\bar{P}_{i,n}x_{i+1}^a T_w H_i^f$ (cf \S\ref{introHbar}),
we obtain finally $M=L$ and we are done.

\medskip
Let us now consider the second isomorphism. Let us fix an adjunction
$(F,E)$ with unit $\eta'$ and counit $\eps'$
and consider the dual categorification $\CA'$ of $\CA(n)$. We denote
by $X'$ and $T'$ its defining endomorphisms.
Define $\sigma':FE\xrightarrow{\eta'FE}EFFE\xrightarrow{ET'E}EFFE
\xrightarrow{EF\eps'}EF$.

Let $G=FE$ and $H=EF$. There is an adjoint pair $(EF,EF)$ with counit
$\eps_H:EFEF\xrightarrow{E\eps'F}EF\xrightarrow{\eps}\Id$ and an
adjoint pair $(FE,FE)$ with unit
$\eta_G:\Id\xrightarrow{\eta} FE\xrightarrow{F\eta'E}FEFE$.
Consider the canonical isomorphism
$$\zeta:
\Hom(FE,EF)=\Hom(G,H)\iso \Hom(H^\vee,G^\vee)\iso \Hom(EF,FE)$$
corresponding to these adjunctions. The commutativity of the
following diagram shows that $\zeta(\sigma')=\sigma$.
$$\xymatrix{
EF\ar[d]_{\eta EF} \\
FEEF\ar[rr]^-{F\eta'EEF}\ar[drr]^-{F\eta'EEF}\ar[rrdd]_-{FEEF} &&
 FEFEEF\ar[rr]^-{FE\eta'FEEF}&& FEEFFEEF\ar[rrd]^-{FEET'EEF} \\
&& FEFEEF\ar[rr]^-{FE\eta'FEEF}\ar[d]^{FE\eps'EF}&&
 FEEFFEEF\ar[rr]^-{FTFFEEF}\ar[d]^{FEEF\eps'EF} && FEEFFEEF
 \ar[d]^{FEEF\eps'EF} \\
&& FEEF \ar[rr]^-{FE\eta'EF} \ar[rrd]_-{FEEF}&& FEEFEF\ar[d]^{FEE\eps'F} &&
 FEEFEF\ar[d]^{FEE\eps'F}\\
&& && FEEF\ar[rr]_-{FTF} && FEEF\ar[d]^{FE\eps}\\
&& && && FE
}$$
Similarly, using the canonical adjoint pair $(\Id,\Id)$, we get
a canonical isomorphism
$$\zeta':\Hom(\Id,EF)=\Hom(\Id,H)\iso \Hom(H^\vee,\Id)
\iso\Hom(EF,\Id).$$
We have
 $\zeta'((\mathbf{1}_E(\tX')^j)\circ \eta')=\eps\circ (X^j \mathbf{1}_F)$.

We have shown that the adjoint to
$$\sigma+\sum_{j=0}^{\lambda-1} \eps\circ(X^j\mathbf{1}_F)\ :\ 
EF\Id_{\CA_{\lambda}}\iso
FE\Id_{\CA_{\lambda}}\oplus \Id_{\CA_{\lambda}}^{\bigoplus\lambda}$$
is 
$$\sigma'+\sum_{j=0}^{\lambda-1} (\mathbf{1}_{F'} (\tX')^j)\circ\eta'\ :\ 
E'F'\Id_{\CA'_{-\lambda}} \oplus \Id_{\CA'_{-\lambda}}^{\bigoplus\lambda}\to
F'E'\Id_{\CA'_{-\lambda}}.$$
One checks easily that the first map of the Proposition remains an
isomorphism if $X_{i+1}$ is replaced by $\tX_{i+1}$.
Since the categorification $\CA'$ is isomorphic to $\CA(n)$, this shows
that the map
$\sigma'+\sum_{j=0}^{\lambda-1} (\mathbf{1}_{F'} (\tX')^j)\circ\eta'$
is an isomorphism, hence 
$\sigma+\sum_{j=0}^{\lambda-1} \eps\circ(X^j\mathbf{1}_F)$ is an isomorphism
as well.
\end{proof}

\subsubsection{}

Let us fix a family $\{M_\lambda\in\CA_\lambda\}_\lambda$. Let
$\CM_\lambda$ be the full subcategory of $\CA_\lambda$ whose objects
are finite direct sums of direct summands of $M_\lambda$.
We assume that
$\CM=\bigoplus_\lambda \CM_\lambda$ is stable under $E$ and $F$.

Let $A'_\lambda=\End_{\CA}(M_\lambda)$,
$\CA'_\lambda=A'_\lambda\mMod$ and
$\CA'=\bigoplus_\lambda \CA'_\lambda$.
We put
\begin{align*}
E'&=\bigoplus_\lambda \Hom_{\CA}(M_{\lambda+2},EM_\lambda)
\otimes_{A'_\lambda}-:\CA'\to\CA'\\
\text{ and }
F'&=\bigoplus_\lambda \Hom_{\CA}(M_{\lambda-2},FM_\lambda)
\otimes_{A'_\lambda}-:\CA'\to\CA'.
\end{align*}

We have $\Hom_{\CA}(M_{\lambda+2},EM_\lambda)\simeq
\Hom_{\CA}(FM_{\lambda+2},M_\lambda)$ and $FM_{\lambda+2}\in\CM_\lambda$.
It follows that
$\Hom_{\CA}(M_{\lambda+2},EM_\lambda)$ is a projective right
$A'_\lambda$-module, so $E'$ is an exact functor. Similarly, $F'$ is an
exact functor. Also, they send projectives to projectives.

\smallskip
Consider the functor
$R=\bigoplus_\lambda M_\lambda\otimes_{A'_\lambda}-:
\CA'\to \CA$. Its restriction to $\CA'\mproj$ is
an equivalence
$\CA'\mproj\iso \CM$. So,
the functor $G\mapsto RG$ from the category of 
exact functors $\CA'\to\CA'$ sending projectives to
projectives to the category of functors $\CA'\to\CA$ is fully faithful.

\smallskip
The canonical map
$$M_{\lambda+2}\otimes_{A'_{\lambda+2}}\Hom_\CA(M_{\lambda+2},EM_\lambda)
\iso EM_\lambda, \ m\otimes f\mapsto f(m)$$
is an isomorphism, since $E M_\lambda\in\CM_{\lambda+2}$.
The induced map
$$M_{\lambda+2}\otimes_{A'_{\lambda+2}}\Hom_\CA(M_{\lambda+2},EM_\lambda)
\otimes_{A'_\lambda}U\iso
E(M_\lambda\otimes_{A'_\lambda}U), \
m\otimes f\otimes u\mapsto E(m'\mapsto m'\otimes u)(f(m))$$
for $U\in A'_\lambda\mMod$ is an isomorphism, since it is
an isomorphism for $U=A'_\lambda$.

We obtain an isomorphism $RE'\iso ER$ and we construct similarly
an isomorphism $RF'\iso FR$.

\smallskip
Let $X'$ (resp. $T'$)
be the inverse image of $X\id_R$ (resp. $T\id_R$) via the canonical
isomorphisms $\End(E')\iso \End(RE')\iso \End(ER)$ (resp.
$\End(E^{'2})\iso \End(RE^{'2})\iso \End(ERE')\iso \End(E^2R)$).

Proceeding similarly, the adjoint pair $(E,F)$ gives an adjoint pair
$(E',F')$ and the functor $F'$ is isomorphic to a left adjoint of $E'$.

\begin{thm}
\label{restrictioncat}
The data above defines an $\Gsl_2$-categorification on $\CA'$ and a morphism
of $\Gsl_2$-categorifications $\CA'\to\CA$.
\end{thm}

\begin{proof}
The $\Gsl_2$-relations in $K_0(\CA'\mproj)$ hold thanks to
Theorem \ref{decH} applied to the restriction of functors to
$\CM$. The local finiteness follows from the case of $\CA$.
The commutativity of the diagrams of Lemma \ref{morphismweak} follows
immediately
from the construction of the adjoint pair $(E',F')$. This shows that
$\CA'$ is a weak categorification and that $R$ defines a morphism of
weak categorifications.

By construction, this weak categorification is a categorification and
the morphism of weak categorifications is actually a morphism of
categorifications.
\end{proof}

\begin{cor}
\label{subcat}
Let $M\in\CA$. Then, there exists a finite dimensional algebra $A$,
an $\Gsl_2$-categorification on $A\mMod$ and a morphism
of $\Gsl_2$-categorifications $R:A\mMod\to\CA$ such that
$M$ is a direct summand of $R(A)$.
\end{cor}

\begin{proof}
Let $N=\bigoplus_{i,j\ge 0}E^iF^j M$, a finite sum. Let $N_\lambda$ be the
projection of $N$ on $\CA_\lambda$. Now, we can apply the constructions
and results above, the stability being provided by Corollary \ref{homotopyK0}.
\end{proof}

\section{Categorification of the reflection}
\label{secrefl}
\subsection{Rickard's complexes}
\label{Rickardcplx}

Let $\lambda\in\BZ$. We construct a complex of
functors 
$$\Theta_\lambda:\Comp(\CA_{-\lambda})\to\Comp(\CA_{\lambda}),$$
following Rickard \cite{Ri1} (originally, for blocks
of symmetric groups).

We denote by $(\Theta_\lambda)^{-r}$ the restriction of
$E^{(\sgn,\lambda+r)}F^{(1,r)}$ to $\CA_{-\lambda}$
for $r,\lambda+r\ge 0$ and we put $(\Theta_\lambda)^{-r}=0$ otherwise.

Consider the map
$$f:E^{\lambda+r}F^r=
E^{\lambda+r-1}EFF^{r-1}\xrightarrow{\mathbf{1}_{E^{\lambda+r-1}}\eps
\mathbf{1}_{F^{r-1}}}E^{\lambda+r-1}F^{r-1}.$$
We have $E^{(\sgn,\lambda+r)}=E^{\lambda+r}
c^{\sgn}_{[\GS_{\lambda+r}/ \GS_{[2,\lambda+r]}]} c_{[2,\lambda+r]}^{\sgn}
\subseteq E^{(\sgn,\lambda+r-1)}E$ and similarly
$F^{(1,r)}\subseteq FF^{(1,r-1)}$, hence $f$ restricts to a map
$$d^{-r}:E^{(\sgn,\lambda+r)}F^{(1,r)}\to
E^{(\sgn,\lambda+r-1)}F^{(1,r-1)}.$$

We put
$$\Theta_\lambda=\cdots\to (\Theta_\lambda)^{-i}\xrightarrow{d^{-i}}
(\Theta_\lambda)^{-i+1}\to\cdots.$$

\begin{lemma}
\label{isacomplex}
$\Theta_\lambda$ is a
complex.
The map
$[\Theta_\lambda]:V_{-\lambda}=K_0(\CA_{-\lambda})\to
V_{\lambda}=K_0(\CA_{\lambda})$ 
coincides with the action of $s$.
\end{lemma}

\begin{proof}
The map $d^{1-r}d^{-r}$ is the restriction of
$\mathbf{1}_{E^{\lambda+r-2}}\eps_2\mathbf{1}_{F^{r-2}}$, where
$\eps_2:EEFF\xrightarrow{\mathbf{1}_{E}\eps\mathbf{1}_F}EF
\xrightarrow{\eps}\Id$.
Since $c^{\sgn}_{\lambda+r}=c^{\sgn}_{[\GS_{\lambda+r}/\GS_2]}c^{\sgn}_2$
and $c^1_r=c_2^1 c^1_{[\GS_2\setminus\GS_r]}$, it follows that
$$E^{(\sgn,\lambda+r)}F^{(1,r)}\subseteq E^{\lambda+r-2}E^{(\sgn,2)}F^{(1,2)}
F^{r-2}.$$
So, in order to prove that $d^{1-r}d^{-r}=0$, it is enough
to show that the composition
$$E^2F^2\xrightarrow{c^{\sgn}_2c^1_2}E^2F^2\Rarr{\eps_2}\Id$$
 vanishes,
where $c^{\sgn}_2$ acts on $E^2$ and $c^1_2$ acts on $F^2$.
This composition is equal to the composition
$E^2F^2\xrightarrow{(c^{\sgn}_2c^1_2)\mathbf{1}_{F^2}}E^2F^2\Rarr{\eps_2}\Id$,
where $c^{\sgn}_2c^1_2$ acts now on $E^2$. We are done, since
$c^{\sgn}_2c^1_2=0$.

The last statement is given by Lemma \ref{simplerrefl}.
\end{proof}

\begin{rem}
Let $M\in\CA_{-\lambda}$.
Let $l=\max\{r\ge 0| F^r M\not=0\}$, a finite integer.
Then, $(\Theta_\lambda)^{-i}(M)=0$ when
$i\not\in [\max(0,-\lambda),l]$.
\end{rem}

\subsection{Derived equivalence from the simple reflection}
\label{dersimple}
Let $\Theta=\bigoplus_\lambda \Theta_{\lambda}$.

The following Lemma follows easily from Lemma \ref{morphismweak}.

\begin{lemma}
\label{morphismreflection}
Let $R:\CA'\to\CA$ be a morphism of $\Gsl_2$-categorifications.
Then, there is an isomorphism of complexes of functors
$\Theta R\iso R\Theta'$.
\end{lemma}

We can now state our main Theorem

\begin{thm}
\label{mainthm}
The complex of functors $\Theta$ induces a self-equivalence
of $K^b(\CA)$ and of $D^b(\CA)$ and
induces by restriction equivalences
$K^b(\CA_{-\lambda})\iso K^b(\CA_\lambda)$ and
$D^b(\CA_{-\lambda})\iso D^b(\CA_\lambda)$.
Furthermore, $[\Theta]=s$.
\end{thm}

\begin{rem}
In the context of symmetric groups, the invertibility of
$\Theta_\lambda$ when the complex has only one (resp. two)
non-zero term 
is due to Scopes \cite{Sco} (resp. Rickard \cite{Ri1}).
\end{rem}

\begin{proof}[Proof of Theorem \ref{mainthm}]
Since $E$ and $F$ have right adjoints,
there is a complex of functors $\Theta_\lambda^\vee$ 
that gives a right adjoint to $\Theta_\lambda$ (cf \S\ref{adjcplx}).
Let $\eps:\Theta_\lambda\Theta_{\lambda}^\vee\to \Id$ be the counit of
adjunction and $Z$ its cone. So, $Z$ is a complex of
exact functors $\CA_{-\lambda}\to\CA_{\lambda}$.

Pick $U\in\CA$ with $FU=0$ and $E^iU\in\CA_{-\lambda}$ and put
$n=h_+(U)$.
The fully faithful functor 
$R_U:K^b(\CA(n)\mproj)\to K^b(\CA)$
commutes with $\Theta_\lambda$ (Lemma \ref{morphismreflection}),
hence commutes with $\Theta_\lambda^\vee$
and with $Z$ (cf \S \ref{commutadj}). By Theorem
\ref{derminimal}, we have $Z(E^iU)=0$. Now,
Proposition \ref{vanishingfunctor} shows that
$Z(M)=0$ in $D^b(\CA_{-\lambda})$ for all $M\in D^b(\CA_{-\lambda})$.
So, $\eps$ is an isomorphism in $D(\CA_{-\lambda})$.
One shows similarly that $\Theta_\lambda$ has a left inverse in
$D^b(\CA_{-\lambda})$.

\smallskip
Let us now prove that $\eps$ is still an isomorphism in $K^b(\CA_{-\lambda})$.
Let $M\in \Comp^b(\CA_{-\lambda})$. By Corollary \ref{subcat},
there is a finite dimensional $k$-algebra $A$,
an $\Gsl_2$-categorification on $\CA'=A\mMod$ and a morphism
of $\Gsl_2$-categorifications $R:\CA'\to\CA$ such that 
the terms of $M$ are direct summands of $R(A)$.
The functor $R$ induces a fully faithful
triangulated functor $R':K^b(\CA'\mproj)\to K^b(\CA)$.
The derived category case of the Theorem shows that $\eps'$ is an
isomorphism in
$K^b(\CA'_{-\lambda}\mproj)\iso D^b(\CA'_{-\lambda})$.
As above, we deduce that
$\eps$ is an isomorphism in the image of $R'$, hence $\eps(M)$
is an isomorphism in $K^b(\CA_{-\lambda})$. One proceeds similarly
to show that $\Theta_\lambda$ has a left inverse in
$K^b(\CA_{-\lambda})$.
\end{proof}

\subsection{Equivalences for the minimal categorification}
\label{equivmin}

\begin{thm}
\label{derminimal}
Let $n\ge 0$ and $\CA=\CA(n)$ be the minimal categorification.
Fix $\lambda\ge 0$ and let $l=\frac{n-\lambda}{2}$.
The homology of the complex of functors $\Theta_\lambda$
is concentrated in degree $-l$ and
$H^{-l}\Theta_\lambda:\CA_{-\lambda}\iso\CA_{\lambda}$
is an equivalence.
\end{thm}

\begin{proof}
In order to show that the homology of $\Theta_\lambda$
is concentrated in degree $-l$, it
suffices to show that $\Theta_\lambda(B_lc_l^1)$ is homotopy
equivalent to a complex concentrated in degree $-l$, since
$B_lc_l^1$ is a progenerator for $B_l\mMod$. This is
equivalent to the property that
$H^*(C)=0$ for $*\not=-l$, 
where $C=c_{n-l}^{\sgn}\bar{H}_{n-l}\otimes_{B_{n-l}}\Theta_\lambda(B_l c_l^1)$,
since $c_{n-l}^{\sgn}\bar{H}_{n-l}$ is the unique simple right $B_{n-l}$-module and
$C^{-r}=0$ for $r>l$.

\smallskip
We have
$$C^{-r}=c_{n-l}^{\sgn}\bar{H}_{n-l}\otimes_{B_{n-l}}B_{n-l}
c_{[l-r+1,n-l]}^{\sgn}
\otimes_{B_{l-r}}c_{[l-r+1,l]}^1 B_l\otimes_{B_l} B_lc_l^1.$$

Lemma \ref{isoB} gives an isomorphism
$$C^{-r}\iso c_{n-l}^{\sgn}\bar{H}_{n-l}
 \otimes_{B_{n-l}}B_{n-l}c_{[l-r+1,n-l]}^{\sgn}
c_{[1,l-r]}^1\otimes_k \bigoplus_{\mu\in P(r,n-l)}
m_\mu(x_{l-r+1},\ldots,x_l)k.$$

Proposition \ref{decomp1} and Lemma \ref{invar} give isomorphisms
$$
\bigoplus_{0\le a_1<\cdots<a_{l-r}<n-l}x_1^{a_1}\cdots x_{l-r}^{a_{l-r}}k
\xrightarrow[\can]{\sim}
\Lambda^{\GS_{l-r}} (P_{n-l}^{\GS_{[l-r+1,n-l]}}
\otimes_{P_{n-l}^{\GS_{n-l}}}k)
$$
\begin{align*}
\Lambda^{\GS_{l-r}} (P_{n-l}^{\GS_{[l-r+1,n-l]}}
\otimes_{P_{n-l}^{\GS_{n-l}}}k)&\iso
c_{n-l}^{\sgn}\bar{H}_{n-l}\otimes_{B_{n-l}}B_{n-l}c_{[l-r+1,n-l]}^{\sgn}
c_{[1,l-r]}^1\\
y&\mapsto c_{n-l}^{\sgn} y
\end{align*}
and these induce isomorphisms
$E^{-r}\xrightarrow[\phi^{-r}]{\sim} D^{-r}\xrightarrow[\psi^{-r}]{\sim}
C^{-r}$, where
$$E^{-r}=\bigoplus_{0\le a_1<\cdots<a_{l-r}<n-l}
x_1^{a_1}\cdots x_{l-r}^{a_{l-r}}k\otimes
\bigoplus_{\mu\in P(r,n-l)} m_\mu(x_{l-r+1},\ldots,x_l)k$$

$$\textrm{and }
D^{-r}=\Lambda^{\GS_{l-r}} (P_{n-l}^{\GS_{[l-r+1,n-l]}}
\otimes_{P_{n-l}^{\GS_{n-l}}}k)
\otimes
\bigoplus_{\mu\in P(r,n-l)} m_\mu(x_{l-r+1},\ldots,x_l)k$$

\smallskip
Let $\mu\in P(r,n-l)$ and $0\le a_1<\cdots<a_{l-r}<n-l$.
Given a positive integer $b$, we write $b\prec \mu$ when $b$
appears in $\mu$ and we denote then by
$\mu\setminus b$ the partition obtained from $\mu$
by removing one instance of $b$.
We have $m_\mu(x_{l-r+1},\ldots,x_l)=
\sum_{b\prec\mu} x_{l-r+1}^bm_{\mu\setminus b}(x_{l-r+2},\ldots,x_l)$.
It follows that
$$d_C^{-r}\psi^{-r}(x_1^{a_1}\cdots x_{l-r}^{a_{l-r}}\otimes m_\mu)=
\psi^{-r+1}\left(\sum_{b\prec\mu}
x_1^{a_1}\cdots x_{l-r}^{a_{l-r}}x_{l-r+1}^b\otimes
m_{\mu\setminus b}\right).$$

Assume $b=n-l$. Since $x_{l-r+1}^{n-l}\in \Gn_{n-l}P_{n-l}$,
it follows that $x_1^{a_1}\cdots x_{l-r}^{a_{l-r}}x_{l-r+1}^b$
is $0$ in
$\Lambda^{\GS_{l-r+1}}
(P_{n-l}^{\GS_{[l-r+2,n-l]}}\otimes_{P_{n-l}^{\GS_{n-l}}}k)$.
One gets the same conclusion when $b\in\{a_1,\ldots,a_{l-r}\}$.

So,
$$d_C^{-r}\psi^{-r}\phi^{-r}(x_1^{a_1}\cdots x_{l-r}^{a_{l-r}}\otimes m_\mu)=
\psi^{-r+1}\phi^{-r+1}\left(
\sum_{b\prec\mu,b\not\in\{a_1,\ldots,a_{l-r},n-l\}}
\sgn(\sigma_b)x_1^{a'_1}\cdots x_{l-r+1}^{a'_{l-r+1}}
\otimes m_{\mu\setminus b}\right)$$
where $\sigma_b\in\GS_{l-r+1}$ is the permutation such that,
putting $a_{l-r+1}=b$ and $a'_j=a_{\sigma_b(j)}$, we have
$a'_1<a'_2<\cdots<a'_{l-r+1}$.

\smallskip
Let $L=k^{n-l}$, with canonical basis $\{e_i\}_{1\le i\le n-l}$.
The Koszul complex $K$ of $L$ is a bigraded $k$-vector space given by
$K^{p,q}=\Lambda^pL\otimes S^qL$, with a differential of bidegree $(-1,1)$
given by
$$(e_{a_1}\cdots e_{a_p})\otimes x\mapsto
\sum_{i=1}^p (-1)^{i+p+1} (e_{a_1}\cdots e_{a_{i-1}}e_{a_{i+1}}\cdots e_{a_p})
\otimes e_{a_i}x.$$

Its dual $\Hom_k(K,k)$ is isomorphic to $J$ defined as follows.
We put $J^{p,q}=\Lambda^p(L^*)\otimes S^q(L^*)$. 
Let $\{f_i\}$ be the dual basis of $L^*$ and 
$f_\mu=f_{\mu(1)}\cdots f_{\mu(q)}\in S^q(L^*)$ for $\mu\in P(q,n-l)$.
Then, the differential 
$d_J:J^{p,q}\to J^{p+1,q-1}$ is given by
$$(f_{a_1}\cdots f_{a_p})\otimes f_\mu\mapsto
\sum_{b\prec\mu,a_1<\cdots<a_i<b<a_{i+1}<\cdots<a_p}
(-1)^{i+p}(f_{a_1}\cdots f_{a_i}f_bf_{a_{i+1}}\cdots
f_{a_p})\otimes f_{\mu\setminus b}.$$

The homology of $J$ is concentrated in bidegree $(0,0)$ and isomorphic
to $k$.
Note that $J^{\bullet,q}$ is a graded right $\Lambda(L^*)$-module, with
action given by right multiplication. This provides $J$ with the
structure of a complex of free graded $\Lambda(L^*)$-modules
(the degree $-q$ term is $J^{\bullet,q}$), hence of
free graded $k[f_{n-l}]/(f_{n-l}^2)$-modules  by restriction. So,
the $(-q)$-th homology group of $J\otimes_{k[f_{n-l}]/(f_{n-l}^2)}k$
is a one-dimensional graded $k$-vector space which is in degree $q$. 
The complexes of vector spaces
$J\otimes_{k[f_{n-l}]/(f_{n-l}^2)}k$ and
$Jf_{n-l}$ are isomorphic, with a shift by one in the grading.
The complex $Jf_{n-l}$ decomposes as the direct sum (over $i$) of the complexes
$\bigoplus_q \Lambda^{i-q}(L^*)f_{n-l}\otimes S^q(L^*)$ and the
cohomology of such a complex is concentrated in degree $-i$.

We have an isomorphism 
$$E^{-r}\iso (\Lambda^{l-r}L^*) f_{n-l}\otimes S^rL^*\subseteq J^{l-r+1,r},\
x_1^{a_1}\cdots x_{l-r}^{a_{l-r}}\otimes m_\mu\mapsto
(f_{a_1}\cdots f_{a_{l-r}}f_{n-l})\otimes f_\mu.$$
This induces an isomorphism between $E$ and the the subcomplex
$\bigoplus_r (\Lambda^{l-r}L^*) f_{n-l}\otimes S^rL^*$ of $J^{l+*+1,-*}$.
It follows that the homology of $E$ is
concentrated in degree $-l$.

\medskip
The complex of functors
$\Theta_{-\lambda}$ is given by tensor product by a bounded complex
of $(B_{n-l},B_l)$-bimodules which are projective as $B_{n-l}$-modules
and as $B_l$-modules. The homology of that complex is concentrated
in the lowest degree where the complex has a non zero component, 
hence the homology $M$ is still projective as a $B_{n-l}$-module
and as a $B_l$-module. Lemma \ref{isacomplex} shows that
$M\otimes_{B_l}-$ sends the unique simple $B_l$-module to the unique
simple $B_{n-l}$-module. By Morita theory,
$M$ induces an equivalence.
\end{proof}

\section{Examples}
\label{secappl}
In this section \S\ref{secappl}, the field $k$ will always be assumed to be big
enough so that the simple modules considered are absolutely simple.

In most of our examples, $\Gsl_2$-categorifications are constructed
in families, using the following recipe.
We start with left and right
adjoint functors $\hat{E}$ and $\hat{F}$
 on an abelian category $\CA$, together with
$X\in\End(\hat{E})$ and $T\in\End(\hat{E}^2)$ satisfying the defining
relations of (possibly degenerate) affine Hecke algebras. We obtain
for each $a\in k$ (with $a\neq 0$ if $q\neq 1$)
an $\Gsl_2$-categorification on $\CA$ given by
$E=E_a$ and $F=F_a$, the generalised $a$-eigenspaces of $X$ acting
on $\hat{E}$ and $\hat{F}$. While we need to check in each example that
$E$ and $F$ do indeed give an action of $\Gsl_2$ on $K_0(\CA)$, it is
automatic that $X$ and $T$ restrict to endomorphisms of $E$ and $E^2$
with the desired properties. That $T$ restricts is a consequence
of the identity (a special case of (1))
\begin{multline*}
T_{1}(X_{2}-a)^{N}-(X_{1}-a)^{N}T_{1} \\
=
\begin{cases}
(q-1)X_2[(X_{1}-a)^{N-1}+(X_{1}-a)^{N-2}
(X_{2}-a)+\cdots+(X_{2}-a)^{N-1}]
& \text{ if }q\not=1 \\
(X_{1}-a)^{N-1}+(X_{1}-a)^{N-2}
(X_{2}-a)+\cdots+(X_{2}-a)^{N-1}
& \text{ if }q=1.
\end{cases}
\end{multline*}
in $H_2(q)$.

\subsection{Symmetric groups}
\label{secsymm}
\subsubsection{}
\label{introsymm}
Let $p$ be a prime number and $k=\BF_p$.
The quotient of $H_n(1)$ by the ideal generated by $X_1$ is the group
algebra $k\GS_n$.
The images of $T_i$ and $X_i$ in $k\GS_n$ are $s_i=(i,i+1)$ and
the Jucys-Murphy element $L_i=(1,i)+(2,i)+\cdots+(i-1,i)$.

Let $a\in k$.
Given $M$ a $k\GS_n$-module, we denote by $F_{a,n}(M)$ the generalized
$a$-eigenspace of $X_n$. This is a $k\GS_{n-1}$-module.
We have a decomposition
$\Res_{k\GS_{n-1}}^{k\GS_n}=\bigoplus_{a\in k}F_{a,n}$.
There is a corresponding decomposition
$\Ind_{k\GS_{n-1}}^{k\GS_n}=\bigoplus_{a\in k}E_{a,n}$, where $E_{a,n}$ is
left and right adjoint to $F_{a,n}$.
We put $E_a=\bigoplus_{n\ge 1} E_{a,n}$ and
$F_a=\bigoplus_{n\ge 1} F_{a,n}$.

\smallskip
Recall the following classical result \cite{LLT}.

\begin{thm}
\label{slp}
The functors $E_a$ and $F_a$ for $a\in \BF_p$ give rise to an action of
the affine Lie algebra
$\widehat{\Gsl}_p$ on $\bigoplus_{n\ge 0}K_0(k\GS_n\mMod)$.

The decomposition of $K_0(k\GS_n\mMod)$ in blocks coincides
with its decomposition in weight spaces.

Two blocks of symmetric groups have the same weight if and only
if they are in the same orbit under the adjoint action of the affine
Weyl group.
\end{thm}

In particular for each $a\in \BF_p$ the functors $E_a$ and $F_a$ give
a weak $\Gsl_2$-categorification on
 $\CA=\bigoplus_{n\ge 0}k\GS_n\mMod$.

We denote by $X$ the endomorphism of $E_a$ given on $E_{a,n}$ by right
multiplication by $L_n$ (on the $(k\GS_n,k\GS_{n-1})$-bimodule
$k\GS_n$).
We denote by $T$ the endomorphism of $E_a^2$ given on $E_{a,n}E_{a,n-1}$
by right multiplication by $s_{n-1}$ (on the $(k\GS_n,k\GS_{n-2})$-bimodule
$k\GS_n$).
This gives an $\Gsl_2$-categorification on $\CA$ (here, $q=1$).

\subsubsection{}

\begin{thm}
\label{equivsymm}
Let $R=k$ or $\BZ_{(p)}$.
Let $A$ and $B$ be two blocks of symmetric groups over $R$ with isomorphic
defect groups. Then, $A$ and $B$ are splendidly Rickard
equivalent (in particular, they are derived equivalent).
\end{thm}

\begin{proof}
Two blocks of symmetric groups over $k$ have isomorphic defect groups
if and only if they have equal weights (cf \S\ref{local} below).
By Theorem \ref{slp}, there is a sequence
of blocks $A_0=A,A_1,\ldots, A_r=B$ such that $A_j$ is the image of
$A_{j-1}$ by some simple reflection $\sigma_{a_j}$ of the affine Weyl group. 
By Theorem \ref{mainthm}, the complex of functors $\Theta$
associated with $a=a_j$ induces a
self-equivalence of $K^b(\CA)$. 
It restricts to a splendid Rickard equivalence between
$A_j$ and $A_{j+1}$. By composing these equivalences, we obtain
a splendid Rickard equivalence between $A$ and $B$
(note that the composition of splendid equivalences can easily be seen
to be splendid, cf eg \cite[Lemma 2.6]{Rou2}).

\smallskip
The constructions of $E$ and $F$ lift uniquely to $\BZ_{(p)}$~:
$\Ind_{\BZ_{(p)}\GS_{n-1}}^{\BZ_{(p)}\GS_n}=\bigoplus_{a\in k}\tilde{E}_a$,
$\Res_{\BZ_{(p)}\GS_{n-1}}^{\BZ_{(p)}\GS_n}=\bigoplus_{a\in k}\tilde{F}_a$,
where $\tilde{E}_a\otimes_{\BZ_{(p)}}k=E_a$,
where $\tilde{F}_a\otimes_{\BZ_{(p)}}k=F_a$ and
$\tilde{E}_a$ and $\tilde{F}_a$ are left and right adjoint.
We denote by $\tilde{T}$ the endomorphism of $\tilde{E}_a^2$ given on
$\tilde{E}_{a,n}\tilde{E}_{a,n-1}$ by the action of $s_{n-1}$.

The construction of $\Theta$ 
in \S \ref{Rickardcplx} lifts to a complex $\tilde{\Theta}$ of functors
on $\tilde{\CA}=\bigoplus_{n\ge 0}\BZ_{(p)}\GS_n\mMod$.
By \cite[end of proof of Theorem 5.2]{Ri3},
the lift $\tilde{\Theta}$ of $\Theta$ is a
splendid self Rickard equivalence of $D^b(\tilde{\CA})$ and we conclude
as before.
\end{proof}

\begin{rem}
The equivalence depends on the choice of a sequence of simple reflections
whose product sends one block to the other. If, as expected,
the categorifications of
the simple reflections give rise to a braid group action on the derived
category of $\bigoplus_{n\ge 0}k\GS_n\mMod$, then one can
choose the canonical lifting of the affine Weyl group element in the
braid group to get a canonical equivalence.
\end{rem}

\begin{rem}
Theorem \ref{equivsymm} gives isomorphisms between Grothendieck groups
of the blocks (taken over $\BQ$) satisfying certain arithmetical properties
(perfect isometries or even isotypies). These arithmetical properties
were shown already by Enguehard \cite{En}.
\end{rem}

\begin{rem}
Two blocks of symmetric groups over $k$ have isomorphic
defect groups if and only if they have the same number of simple
modules, up to the exception of blocks of weights $0$ and $1$ for
$p=2$ --- note that a block of weight $0$ is simple whereas a block of
weight $1$ is not simple, so two such blocks are not derived equivalent.
So, one can restate Theorem \ref{equivsymm} as follows~:

Let $A$ and $B$ be two blocks of symmetric groups over $k$.
Then, $A$ and $B$ are derived equivalent if and only if they have
isomorphic defect groups.
Assume $A$ and $B$ are not simple if $p=2$.
Then, $A$ and $B$ are derived equivalent
if and only if $\rank K_0(A)=\rank K_0(B)$.
\end{rem}

We can now deduce a proof of Brou\'e's abelian defect group conjecture
for blocks of symmetric groups~:

\begin{thm}
\label{mainSn}
Let $A$ be a block of a symmetric group $G$ over $\BZ_{(p)}$,
$D$ a defect group and
$B$ the corresponding block of $N_G(D)$. If $D$ is abelian, then
$A$ and $B$ are splendidly Rickard equivalent.
\end{thm}

\begin{proof}
By \cite{ChKe}, there is a block $A'$ of a symmetric group
which is splendidly Morita equivalent to the principal block of
$\BZ_{(p)}(\GS_p\wr \GS_w)$, where $w$ is the weight of $A$.
We have a splendid Rickard equivalence between
the principal block of $\BZ_{(p)}\GS_p$ and $\BZ_{(p)}N$, where $N$ is
the normalizer of a Sylow $p$-subgroup of $\GS_p$ by \cite[Theorem 1.1]{Rou2}.
By \cite[Theorem 4.3]{Ma} (cf also \cite[Lemma 2.8]{Rou2} for the
Rickard/derived equivalence part),
we deduce a splendid Rickard equivalence between the
principal blocks of $\BZ_{(p)}(\GS_p\wr \GS_w)$ and $\BZ_{(p)}(N\wr \GS_w)$.
Now, we have an isomorphism
$B\simeq \BZ_{(p)}(N\wr \GS_w)\otimes B_0$, where $B_0$ is a matrix algebra
over $\BZ_{(p)}$, hence there is a splendid Morita equivalence between
$B$ and $\BZ_{(p)}(N\wr \GS_w)$.
So, we obtain
a splendid Rickard equivalence between $B$ and $A'$.

By Theorem \ref{equivsymm}, we have a splendid Rickard equivalence between
$A$ and $A'$ and the Theorem follows.
\end{proof}

\begin{rem}
The existence of an isotypy 
between $A$ and $B$ 
in Theorem~\ref{mainSn} was known by \cite{Rou1}.
\end{rem}

\subsubsection{}
\label{local}
Let us analyze more precisely the categorification.

Given $\lambda$ a partition of $m$, we denote by $|\lambda|=m$ the size
of $\lambda$.
Let $\kappa$ be a $p$-core and $n$ an integer such that
$p|(n-|\kappa|)$ and $n\ge |\kappa|$.
We denote by $b_{\kappa,n}$ the corresponding block of $k\GS_n$
(the irreducible characters of that block
are associated to the partitions having $\kappa$ as their $p$-core).
The integer $\frac{n-|\kappa|}{p}$ is the weight of the block (this notion
of weight is not to be confused with the weights relative to Lie algebra
actions).

Let $\lambda$ be a partition with $p$-core $\kappa$ and $\lambda'$
a partition obtained from $\lambda$ by adding an $a$-node. Then, 
the $p$-core of $\lambda'$ depends only on $\kappa$ and $a$ and
we denote it by $e_a(\kappa)$. Similarly, we define $f_a(\kappa)$
by removing an $a$-node.

We will freely identify a functor $M\otimes-$ with the bimodule $M$.
We have 
\begin{equation}
\label{blockE}
E_{a,n+1}=\bigoplus_{\kappa} b_{e_a(\kappa),n+1}k\GS_{n+1}b_{\kappa,n}
\end{equation}
where $\kappa$ runs over the $p$-cores such that
$|\kappa|\le n$, $|\kappa|\equiv n\pmod p$ and $|e_a(\kappa)|\le n+1$.

\smallskip
Let $b_{\kappa_{-r},l},b_{\kappa_{-r+2},l+1},\ldots,b_{\kappa_{r},l+r}$
be a chain of blocks with $|f_a(\kappa_{-r})|>l-1$,
$|e_a(\kappa_r)|>l+r+1$ and $f_a(\kappa_j)=\kappa_{j-2}$.

Put $n_i=l+(i-r)/2$ and
$B_i=k\GS_{n_i}b_{\kappa_{i},n_i}$ for $-r\le i\le r$ and $i\equiv r\pmod 2$.

Let $\CA=\bigoplus_i B_i\mMod$.
The action of $E=E_a$ and $F=F_a$ on $K_0(\CA)$ gives a representation of
$\Gsl_2$.
This gives an $\Gsl_2$-categorification (here, $q=1$).

The complex of functors $\Theta$ 
restricts to a splendid Rickard equivalence between $B_i$ and $B_{-i}$.

\smallskip
Let us recall some results of the local block theory of symmetric groups
(cf \cite{Pu1} or \cite[\S 2]{Br}).

Let $P$ be a $p$-subgroup of $\GS_n$. Up to conjugacy, we can assume
$[1,n]^P=[n_P+1,n]$ for some integer $n_P$ (we call such a $P$ a
standard $p$-subgroup). Then,
$C_{\GS_n}(P)=H\times \GS_{[n_P+1,n]}$ where $H=C_{\GS_{n_P}}(P)$.
The algebra $kH$ has a unique block.

Given $G$ a finite group and $P$ a $p$-subgroup of $G$, we denote
by $\br_P:(kG)^P\to kC_G(P)$ the Brauer morphism (restriction of
the morphism of $k$-vector spaces $kG\to kC_G(P)$ which is the identity
on $C_G(P)$ and $0$ on $G-C_G(P)$).
We denote by $\Br_P:kG\mMod\to kC_G(P)\mMod$ the Brauer functor given
by $M\mapsto M^P/(\sum_{Q<P}\Tr_Q^P M^Q)$, where
$\Tr_Q^P(x)=\sum_{g\in P/Q}g(x)$.

\smallskip
We will use the following result of Puig and Marichal

\begin{thm}
\label{BrauerSn}
We have
$$\br_P(b_{\kappa,n})=\begin{cases}
1\otimes b_{\kappa,n-n_P} & \text{ if } \frac{n-n_P-|\kappa|}{p}\in\BZ_{\ge 0}
\\
0 & \text{ otherwise.}
\end{cases}
$$
\end{thm}

Note in particular that a standard $p$-subgroup $P$ is a defect group of
$b_{\kappa,n}$ if and only if $P$ is a Sylow $p$-subgroup of
$\GS_{n-|\kappa|}$. In particular, two blocks of symmetric groups have
isomorphic defect groups if and only if they have equal weights.

\smallskip
So, we deduce from (\ref{blockE}) and Theorem \ref{BrauerSn}~:

\begin{lemma}
\label{Brauer}
We have an isomorphism of
$((kH\otimes k\GS_{n-n_P+i}),(kH\otimes k\GS_{n-n_P-1}))$-bimodules
$$\Br_{\Delta P}(E_{a,n+i}\cdots E_{a,n+1}E_{a,n})\iso kH\otimes 
E_{a,n-n_P+i}\cdots E_{a,n-n_P+1}E_{a,n-n_P}.$$
For $i=1$, it is compatible with the action of $T$.

Let $P$ be a non-trivial standard $p$-subgroup of $\GS_{n_{-i}}$.
If $\br_P(b_{\kappa_{i},n_i})$ is not $0$, then
$$\Br_{\Delta_P}(b_{\kappa_{-i},n_{-i}} \Theta b_{\kappa_i,n_i}) \simeq
kH\otimes b_{\kappa_{-i},n_{-i}-n_P}\Theta b_{\kappa_i,n_i-n_P}.$$
\end{lemma}

Note that this Lemma permits to deduce a proof of the Rickard
equivalence in Theorem \ref{equivsymm} from that of the derived equivalence,
by induction on the size of the defect group~:
By induction, $b_{\kappa_{-i},n_{-i}-n_P}\Theta b_{\kappa_i,n_i-n_P}$
induces a Rickard equivalence.
Now, $\Theta$ induces a derived equivalence,
so, it follows from Theorem \ref{splendid} below that $\Theta$ induces a
Rickard equivalence between $B_i$ and $B_{-i}$.

\smallskip
If a splendid complex induces local derived equivalences, then it
induces a Rickard equivalence
\cite[Theorem 5.6]{Rou4} (in
a more general version, but whose proof extends with no
modification)~:
\begin{thm}
\label{splendid}
Let $G$ be a finite group, $b$ a block of $kG$ and $D$ a defect group of
$b$. We assume $b$ is of principal type, \ie, $\br_D(b)$ is a block
of $kC_G(D)$. Let $H$ be a subgroup of $G$ containing $D$ and controlling
the fusion of $p$-subgroups of $D$. Let $c$ be the block of $kH$
corresponding to $b$.

Let $C$ be a bounded complex of $(kGb,kHc)$-bimodules. We assume
$C$ is splendid, \ie, the components $M$ of $C$ are direct summands
of modules $\Ind_{\Delta D}^{G\times H^\circ}N$, where $N$ is a
permutation $\Delta D$-module.

Assume
\begin{itemize}
\item
$\Br_{\Delta P}(C)$ induces a Rickard equivalence between
$kC_G(P)\br_P(b)$ and $kC_H(P)\br_P(c)$ for $P$ a non trivial $p$-subgroup of
$D$ and
\item $C$ induces a derived equivalence between $kGb$ and $kHc$.
\end{itemize}

Then, $C$ induces a Rickard equivalence between $kGb$ and $kHc$.
\end{thm}

\subsection{Cyclotomic Hecke algebras}
\label{cyclo}
\subsubsection{}
\label{defcyclo}
We consider here the non-degenerate case $q\not=1$.
We fix $v_1,\ldots,v_d\in k^\times$.

We denote by 
$\CH_n=\CH_n(v,q)$ the quotient of
$H_n(q)$ by the ideal generated by
$(X_1-v_1)\cdots (X_1-v_d)$.
This is the Hecke algebra of the complex
reflection group $G(d,1,n)$ (cf e.g. \cite[\S 13.1]{Ar2}).

The algebra $\CH_n$ is free over $k$ with basis
$\{X_1^{a_1}\cdots X_n^{a_n}T_w\}_{0\le a_i<d, w\in\GS_n}$ \cite{ArKo}. 
In particular $\CH_{n-1}$ embeds as a subalgebra of
$\CH_n$, and $\CH_n$ is free as a left and as
a right $\CH_{n-1}$-module, for the multiplication action.
The algebra $\CH_n$ is symmetric \cite{MalMat}.

\subsubsection{}
\label{catcyclo}
Let $a\in k^\times$. Given $M$ an $\CH_n$-module, we denote by $F_{a,n}M$ the generalized
$a$-eigenspace of $X_n$. This is an $\CH_{n-1}$-module.
We have a decomposition
$\Res_{\CH_{n-1}}^{\CH_n}=\bigoplus_{a\in k^\times}F_{a,n}$.
There is a corresponding decomposition
$\Ind_{\CH_{n-1}}^{\CH_n}=\bigoplus_{a\in k^\times}E_{a,n}$, where $E_{a,n}$ is
left and right adjoint to $F_{a,n}$.
We put $E_a=\bigoplus_{n\ge 1} E_{a,n}$ and
$F_a=\bigoplus_{n\ge 1} F_{a,n}$.

Now fix $a\in k^\times$. The functors $E=E_a$ 
and $F=F_a$ give an action of ${\Gsl}_2$ on
$\bigoplus_{n\ge 0}K_0(\CH_n\mMod)$ in which the classes of
simple modules are weight vectors \cite[Theorem 12.5]{Ar2}
(only the case where each
parameter if a power of $q$ is considered there, but the proof extends
immediately to our more general setting).
We obtain an $\Gsl_2$-categorification on
$\bigoplus_{n\ge 0} \CH_n\mMod$, where
the endomorphism $X$ of $E$ is given on $E_{a,n}$ by right multiplication
by $X_n$, and the endomorphism $T$ of $E^2$ is given on $E_{a,n}E_{a,n-1}$
by right multiplication by $T_{n-1}$.

\begin{rem}
Let $e$ be the multiplicative order of $q$ in $k^\times$.
Fix $a_0\in k^\times$ and let
$I=\{q^m a_0 \mid m\in\BZ\}$. Then the functors $E_a$ and $F_a$ for $a\in I$ define
an action of $\widehat{\Gsl}_e$ on $\bigoplus_{n\ge 0}K_0(\CH_n\mMod)$.
\end{rem}

\subsubsection{}
Consider here the case $d=1$. Then, $\CH_n=\CH_n(1,q)$ is the Hecke algebra
of $\GS_n$.
Let $e$ be the multiplicative order of $q$ in $k$.
We have a notion of weight of a block as in \S\ref{introsymm}, replacing
$p$ by $e$ in the definitions.

We obtain a $q$-analog of Theorem \ref{equivsymm}:
\begin{thm}
\label{equivHecke}
Assume $d=1$. Let $A$ be a block of $\CH_n$ and $B$ a block of
$\CH_m$. Then, $A$ and $B$ are derived equivalent if and only
they are Rickard equivalent if and only
if they have the same weight.
\end{thm}

\begin{rem}
All of the constructions and results of \S\ref{cyclo} hold for degenerate
cyclotomic Hecke algebras as well, under the assumption that they are symmetric
algebras (which should be provable along the lines of \cite{MalMat}).
\end{rem}

\subsection{General linear groups over a finite field}
\label{GLq}
\subsubsection{}
\label{GLqcat}
Let $q$ be a prime power, $n\ge 0$ and $G_n=\GL_n(q)$. We assume
that $k$ has characteristic $\ell>0$ and $\ell{\!\not|} q(q-1)$.
Let $A_n=kG_nb_n$ be the sum of the unipotent blocks of $kG_n$.

Given $H$ a finite group, we put $e_H=\frac{1}{|H|}\sum_{h\in H}h$.
We denote by ${^tg}$ the transpose of a matrix $g$.

We denote by $V_n$ the subgroup of upper triangular
matrices of $G_n$ with diagonal coefficients $1$ and whose off-diagonal
coefficients vanish outside the $n$-th column. We denote by $D_n$ the subgroup
of $G_n$ of diagonal matrices with diagonal entries $1$ except the $(n,n)$-th one.
$$V_n=\begin{pmatrix}
1 &      & & *\\
  &\ddots& & \vdots \\
  &      &1& * \\
  &      & & 1 \\
\end{pmatrix},\ \
D_n=\begin{pmatrix}
1 \\
 &\ddots\\
 &      &1\\
 &      & &*\\
\end{pmatrix}
$$
Let $i\le n$. We view $G_i$ as a subgroup of $G_n$ via the first $i$
coordinates.

\smallskip
We put
\begin{align*}
E_{i,n}&=kG_ne_{(V_n\rtimes\cdots\rtimes V_{i+1})\rtimes (D_{i+1}\times\cdots
\times D_n)}\otimes_{kG_i}-:A_i\mMod\to A_n\mMod\\
\textrm{ and }
F_{i,n}&=e_{(V_n\rtimes\cdots\rtimes V_{i+1})\rtimes (D_{i+1}\times\cdots\times
D_n)}kG_n\otimes_{kG_n}-:A_n\mMod\to A_i\mMod.
\end{align*}
These functors are canonically left
and right adjoint. Furthermore, there are canonical
isomorphisms $E_{j,n}\circ E_{i,j}\iso E_{i,n}$ and
$F_{i,j}\circ F_{j,n}\iso F_{i,n}$ for $i\le j\le n$.

Let $\CA=\bigoplus_{n\ge 0}A_n\mMod$,
$E=\bigoplus_{n\ge 0}E_{n,n+1}$ and 
$F=\bigoplus_{n\ge 0}F_{n,n+1}$.

\smallskip
We denote by $T$ the endomorphism of $E$ given on
$E_{n-2,n}$ by right multiplication by
$$\hat{T}_{n-1}=q e_{V_nV_{n-1}D_{n-1}D_n} (n-1,n)e_{V_nV_{n-1}D_{n-1}D_n}.$$
We denote by $X$ the endomorphism of $E^2$ given on
$E_{n-1,n}$ by right multiplication by 
$$\hat{X}_n=q^{n-1}e_{V_nD_n}e_{^tV_n}e_{V_nD_n}.$$

\begin{lemma}
\label{GLH}
We have 
$$
(\mathbf{1}_ET)\circ(T\mathbf{1}_E)\circ(\mathbf{1}_ET)=
(T\mathbf{1}_E)\circ(\mathbf{1}_ET)\circ(T\mathbf{1}_E),\quad
(T+\mathbf{1}_{E^2})\circ(T-q\mathbf{1}_{E^2})=0$$
$$\textrm{ and }\quad T\circ(\mathbf{1}_E X)\circ T=
qX\mathbf{1}_E.$$ 
\end{lemma}

\begin{proof}
The first statements involving only $T$'s are the classical results
of Iwahori.

Let $U$ be the subgroup of $G_n$ with diagonal coefficients $1$ and
whose off-diagonal coefficients vanish except the $(n,n-1)$-th.
We have
$$\hat{T}_{n-1}\hat{X}_{n-1}\hat{T}_{n-1}=
q^{n}e_{V_nV_{n-1}D_{n-1}D_n}e_U(n-1,n)e_{^tV_{n-1}}
(n-1,n)e_Ue_{V_{n-1}V_nD_{n-1}D_n}$$
$$=q^{n}e_{V_nV_{n-1}D_{n-1}D_n}e_{^tV_n}e_{V_nV_{n-1}D_{n-1}D_n}=
qe_{V_{n-1}D_{n-1}}\hat{X}_ne_{V_{n-1}D_{n-1}}$$
and this induces the same endomorphism of $E_{n-2,n}$ as $q\hat{X}_n$.
\end{proof}

Lemma \ref{GLH} shows that we have a morphism
$H_n(q)\to \End(E_{0,n})=\End_{kG_n}(kG_n/B_n)$
which sends $T_i$ to the endomorphism given by right multiplication
by $qe_{B_n}(i-1,i)e_{B_n}$ and $X_1$ to the identity,
where $B_n$ is the subgroup of $G_n$ of upper
triangular matrices (cf \S~\ref{actiononpowers}). The classical
result of Iwahori states that the restriction
of this morphism to $H_n^f$ is an isomorphism. This gives us a surjective
morphism $p:H_n\to H_n^f$ whose restriction to $H_n^f$ is the identity.
Since $X_1$ maps to $1$ in $\End(E_{0,n})$ and
the quotient of $H_n$ by $X_1-1$ is isomorphic to $H_n^f$, it follows that $p$
is the canonical map $H_n\to H_n^f$. In particular, the image of $X_i$
is (up to an affine transformation) a Jucys-Murphy element~:
$$p(X_i)=q^{1-i}T_{i-1}\cdots T_1 T_1\cdots T_{i-1}=
1+q^{1-i}(q-1)\left(T_{(1,i)}+qT_{(2,i)}+\cdots+q^{i-2}T_{(i-1,i)}\right).$$

\smallskip
We put $R_n=\Hom_{kG_n}(kG_ne_{B_n},-)
=e_{B_n}kG_n\otimes_{kG_n}-:A_n\mMod\to H_n^f\mMod$.
The multiplication maps
$$e_{B_i}kG_i\otimes_{kG_i}e_{V_n\cdots V_{i+1}D_n\cdots D_{i+1}}kG_n
\to e_{B_n}kG_n$$
and
$$e_{B_n}kG_ne_{B_n}\otimes_{e_{B_i}kG_ie_{B_i}}e_{B_i}kG_i\to
e_{B_n}kG_ne_{V_n\cdots V_{i+1}D_n\cdots D_{i+1}}$$
are isomorphisms. They induce isomorphisms of functors
$$R_i F_{i,n}\iso \Res_{H_i^f}^{H_n^f} R_n \qquad\textrm{and}\qquad
\Ind_{H_i^f}^{H_n^f} R_i\iso R_n E_{i,n}.$$

\begin{rem}
The constructions carried out here make sense more generally for finite groups
with a BN-pair and for arbitrary standard parabolic subgroups, the transpose
operation corresponding to passing from the unipotent radical of a parabolic
subgroup to the unipotent radical of the opposite parabolic subgroup. This
produces a very general kind of ``Jucys-Murphy element'' in Hecke algebras of
finite Weyl groups. In type $B$ or $C$, we should recover the usual
Jucys-Murphy elements.
\end{rem}

Given $a\in k^\times$, let $E_a$ be the generalized $a$-eigenspace of $X$ acting on
$E$.

\begin{lemma}
\label{sl2relGL}
The action of $[E_a]$ and $[F_a]$ on $\bigoplus_{n\ge 0}K_0(A_n\mMod)$ gives
a representation of $\Gsl_2$. Furthermore, the classes of simple objects
are weight vectors.
\end{lemma}

\begin{proof}
Let $\CO$ be a complete discrete valuation ring with field of fractions $K$ and
residue field $k$.
We consider the setting above where $k$ is replaced by $K$.
The functor
$\Hom_{KG_n}(KG_n e_{B_n},-)$ induces an isomorphism from the Grothendieck
group $L_n$ of the category of unipotent representations of $KG_n$ to the
Grothendieck group of the category of representations of the Hecke algebra
of type $\GS_n$ with parameter $q$ over $K$. This isomorphism is compatible
with the actions of $E_a$ and $F_a$. It follows from \S \ref{catcyclo} that $E_a$
and $F_a$ give a representation of $\Gsl_2$ on $\bigoplus_{n\ge 0}L_n$ and
the class of a simple unipotent representation of $KG_n$ is a weight vector.
Now, the decomposition map $L_n\to K_0(A_n)$ is an isomorphism
\cite[Theorem 16.7]{Jam} and the
result follows.
\end{proof}

So, we have constructed an $\Gsl_2$-categorification on $\bigoplus_{n\ge 0}
A_n\mMod$ and a morphism of $\Gsl_2$-categorifications 
$\bigoplus_{n\ge 0} A_n\mMod\to \bigoplus_{n\ge 0}H_n^f\mMod$.

\begin{rem}
Note that we deduce from this that the blocks of $A_n$ correspond to
the blocks of $H_n^f$.
\end{rem}

\subsubsection{}
We assume here only that $\ell{\!\not|}q$.
Let $\CO$ be the ring of integers of a finite extension of $\BQ_\ell$ and
$k$ be the residue field of $\CO$.

Let us recall \cite{FoSri} that the $\ell$-blocks of $\GL_n(q)$
are parametrized by pairs $((s),(B_1,\ldots,B_r))$ where $s$ is a conjugacy
class of semi-simple $\ell'$-elements of $\GL_n(q)$ and
$B_i$ is a block of $\CH_{n_i}(q^{d_i})$,
where $C_{\GL_n(q)}(s)=\GL_{n_1}(q^{d_1})\times\cdots\times
\GL_{n_r}(q^{d_r})$. Let $w_i$ be the $e_i$-weight of
the block $B_i$, where $e_i$ is the multiplicative order of $q^{d_i}$ in
$k^\times$. We define the weight of the block as the family
$\{(w_i,d_i)\}_{1\le i\le r}$.

\begin{thm}
Let $R=k$ or $\CO$.
Two blocks of general linear groups over $R$ with same weights
are splendidly Rickard equivalent.
\end{thm}

\begin{proof}
The results on the local block theory of symmetric groups generalize
to unipotent blocks of general linear groups \cite[\S 3]{Br} and
we conclude as in the proof of Theorem \ref{equivsymm} that the
Theorem holds for unipotent blocks.

By \cite{BoRou2}, a block of a general linear group is
splendidly Rickard equivalent to a unipotent block of 
a product $\GL_{n_1}(q^{d_1})\times\cdots\times \GL_{n_r}(q^{d_r})$
(\cite[Th\'eor\`eme B]{BoRou1} already
provides a complex with homology only in one degree inducing a Morita
equivalence).
Such a block is splendidly Rickard equivalent to the principal
block of
$\GL_{e_1w_1}(q^{d_1})\times\cdots\times \GL_{e_rw_r}(q^{d_r})$ by the
unipotent case of the Theorem.
\end{proof}

\begin{rem}
Assume $l|(q-1)$. Then, $k\GL_n(q)$ has a unique unipotent block, the principal
block. The number of simple modules for such a block is the number of
partitions of $n$. Consequently, a unipotent block
of $\GL_n(q)$ is not derived equivalent to a unipotent block of
$\GL_m(q)$ when $n\not=m$.
\end{rem}

\begin{thm}
Let $A$ be a block of a general linear group $G$ over $R=k$ or $\CO$, let
$D$ be a defect group and
$B$ the corresponding block of $N_G(D)$. If $D$ is abelian, then
$A$ and $B$ are splendidly Rickard equivalent.
\end{thm}

\begin{proof} By the result of \cite{BoRou2} stated above, we may
assume that $A$ is a unipotent block. Then
we proceed as in the proof of Theorem \ref{mainSn}, using
the fact that there is a unipotent block of a general linear group
with defect group isomorphic to $D$
that is splendidly Morita equivalent to the principal block of
$R(\GL_e(q)\wr \GS_w)$ for some $w\ge 0$, where $e$ is the
order of $q$ in $k^\times$ \cite{Pu2, Mi, Tu}.
\end{proof}

\subsection{Category $\CO$} \label{CatO}

\subsubsection{}
We construct here $\Gsl_2$-categorifications on category
$\CO$ of $\Ggl_n$. In particular we show that the
weak $\Gsl_2$-categorification on singular blocks given by Bernstein,
Frenkel and Khovanov \cite{BeFreKho}
is an $\Gsl_2$-categorification.

We denote by $\Gh$ the Cartan subalgebra of diagonal matrices and
$\Gn$ the nilpotent algebra of strictly upper trangular matrices of
the complex Lie algebra $\Gg=\Ggl_n$. We denote by
$\CO$ the BGG category of finitely generated $U(\Gg)$-modules
that are diagonalisable for $\Gh$ and locally nilpotent for
$U(\Gn)$. 

Let $\{e_{ij}\}$ be the standard basis of $\Gg$, and let
$\varepsilon_{1},\ldots,\varepsilon_{n}$ be the basis of $\Gh^{\ast}$
dual to $e_{11},\ldots,e_{nn}$. For each $\lambda\in\Gh^{\ast}$ we
denote by $\lambda_{1},\ldots,\lambda_{n}$ the coefficients of $\lambda$ with
respect to $\varepsilon_{1},\ldots,\varepsilon_{n}$. We write 
$\lambda\to_{a}\mu$ if there exists $j$ such that $\lambda_{j}-j+1=a-1$,
$\mu_{j}-j+1=a$ and $\lambda_{i}=\mu_{i}$ for $i\neq j$.

For each $\lambda\in \Gh^{\ast}$ let $M(\lambda)$ be the Verma module with
highest weight $\lambda$ and let $L(\lambda)$ be its unique irreducible
quotient. Recall that
$M(\lambda)=U(\Gg)\otimes_{U(\Gb)}\BC_{\lambda}$, where $\Gb$ is the subalgebra of
upper-triangular matrices and $\BC_{\lambda}$ is the one-dimensional
$\Gb$-module on which $e_{ii}$ acts as multiplication by 
$\lambda_{i}$. 

Let $\Theta$ be the set of maximal ideals of the center $Z$ of $U(\Gg)$.
For each $\theta\in\Theta$ denote by $\CO_{\theta}$ the full
subcategory of $\CO$ consisting of modules annihilated by some power
of $\theta$. The category $\CO$ splits as a direct sum of the
subcategories $\CO_{\theta}$. Let $\pr_\theta:\CO\to\CO$
denote the projection onto $\CO_\theta$. Each Verma module belongs to some
$\CO_{\theta}$, and $M(\lambda)$ and $M(\mu)$ belong to the same
subcategory if and only if
$\lambda$ and $\mu$ are in the same orbit 
in the dot action of the Weyl group of $\Gg$ on $\Gh^*$, \ie, 
if and only if
 $(\lambda_{1},\lambda_{2}-1,\ldots,\lambda
_{n}-n+1)$ and $(\mu_{1},\mu_{2}-1,\ldots,\mu_{n}-n+1)$ are in the same
$\mathfrak{S}_{n}$-orbit.
 We write $\theta\to_{a}\theta^{\prime}$ if
there exist $\lambda,\mu\in\Gh^{\ast}$ such that 
$M(\lambda)\in\CO_{\theta}$,
$M(\mu)\in\CO_{\theta^{\prime}}$ and
$\lambda\to_{a}\mu$. 

Let $V$ be the natural $n$-dimensional representation of $\Gg$. The
functor $V\otimes -:\CO\to\CO$ decomposes as a direct sum
$\bigoplus_{a\in\BC}E_{a}$, where
$$E_{a}=
\bigoplus\limits_{\substack{\theta,\theta^{\prime}\in
\Theta\\\theta\to_{a}\theta^{\prime}}}
\pr_{\theta^{\prime}}\circ(V\otimes -)\circ \pr_{\theta}.
$$
Each summand $E_{a}$ has a left and right adjoint
$$F_{a}=
\bigoplus\limits_{\substack{\theta,\theta^{\prime}\in
\Theta\\\theta\to_{a}\theta^{\prime}}}
\pr_{\theta}\circ(V^{\ast}\otimes -)\circ \pr_{\theta^{\prime}}.
$$

Let $\lambda\in\Gh^{\ast}$. We have $V\otimes
M(\lambda)=V\otimes(U(\Gg)\otimes_{U(\Gb)}\BC_{\lambda})
\simeq U(\Gg)\otimes_{U(\Gb)}(V\otimes\BC_{\lambda})$, and
therefore $V\otimes M(\lambda)$ has a 
filtration with quotients isomorphic to
the modules $M(\lambda+\varepsilon_{i})$, $i=1,\ldots,n$.  Similarly
$V^{\ast}\otimes M(\lambda)$ has a filtration
with quotients isomorphic to the modules
$M(\lambda-\varepsilon_{i})$, $i=1,\ldots,n$. It follows that
$$
[E_{a}M(\lambda)]=\sum_{\substack{\mu\in\Gh^{\ast}\\
\lambda\to_{a}\mu}}[M(\mu)],\qquad [F_{a}M(\lambda
)]=\sum_{\substack{\mu\in\Gh^{\ast}\\
\mu\to_{a}\lambda}}[M(\mu)]
$$
in $K_{0}(\CO)$. Hence
$$
[E_{a}F_{a}M(\lambda)]-[F_{a}E_{a}M(\lambda)]=c_{\lambda,a}
[M(\lambda)],
$$
where
$c_{\lambda,a}=\#\{i\mid\lambda_{i}-i+1=a\}-\#\{i\mid\lambda_{i}-i+1=a-1\}$.
Because the classes of Verma modules are a basis for $K_{0}(\mathcal{O)}$, we
deduce that for each $a\in\BC$ the functors $E_{a}$ and $F_{a}$
give a weak $\Gsl_2$-categorification on $\CO$ in which
the simple module $L(\lambda)$ has
weight $c_{\lambda,a}$. 

\subsubsection{}
Given $M$ a $\Gg$-module, we have an action map
$\Gg\otimes M \to M$. Let $X_M\in\End_\Gg(V\otimes M)$ be the corresponding
adjoint map. This defines an endomorphism $X$ of the functor $V\otimes -$.
We have $X_M(v\otimes m)=\Omega(v\otimes m)$
where $\Omega=\sum_{i,j=1}^{n}e_{ij}\otimes e_{ji}\in\Gg
\otimes\Gg$.

Define $T_{M}\in\End_{\Gg}(V\otimes V\otimes M)$
by $T_{M}(v\otimes v^{\prime}\otimes m)=v^{\prime}\otimes v\otimes m$.
This defines an endomorphism $T$ of the functor $V\otimes V\otimes -$.

\begin{lemma}
\label{Orelation}
We have the following equality in $\End_\Gg
(V\otimes V\otimes M)$:
$$
T_{M}\circ(\mathbf{1}_{V}\otimes X_{M})=X_{V\otimes M}\circ T_{M}
-\mathbf{1}_{V\otimes V\otimes M}.
$$
\end{lemma}

\begin{proof}
We have
\begin{align*}
X_{V\otimes M}T_{M}(v\otimes v^{\prime}\otimes m)  & =\sum_{i,j=1}^{n}
e_{ij}v^{\prime}\otimes e_{ji}(v\otimes m)\\
& =\sum_{i,j=1}^{n}e_{ij}v^{\prime}\otimes e_{ji}v\otimes m+\sum_{i,j=1}
^{n}e_{ij}v^{\prime}\otimes v\otimes e_{ji}m\\
& =v\otimes v^{\prime}\otimes m+T_{M}(\mathbf{1}_{V}\otimes X_{M})(v\otimes
v^{\prime}\otimes m).
\end{align*}
\end{proof}

The lemma implies that for each $l$ we can define a morphism $H_{l}
(1)\to\End_{\Gg}(V^{\otimes l}\otimes M)$
by $T_{i}\mapsto\mathbf{1}_{V}^{\otimes l-i-1}\otimes T_{V^{\otimes
i-1}\otimes M}$ and $X_{i}\mapsto\mathbf{1}_{V}^{\otimes l-i}\otimes
X_{V^{\otimes i-1}\otimes M}$. Jon Brundan has pointed out to us that this
coincides up to shift (cf Remark~\ref{scalarshift})
with an action described by Arakawa and Suzuki \cite[\S2.2]{ArSu}.

\subsubsection{}
We shall now show that $X$ and $T$ restrict to give endomorphisms of the
functors $E_{a}$ and $E_{a}^{2}$ which define $\Gsl_2$-categorifications
on $\CO$. In view of Lemma~\ref{Orelation}, it suffices to identify
$E_a$ as the generalised $a$-eigenspace of $X$ acting on $V\otimes -$.

To this end we observe that
$
\Omega=\frac{1}{2}(\delta(C)-C\otimes 1 - 1 \otimes C),
$
where $C=\sum_{i,j=1}^{n}e_{ij}e_{ji}\in Z$ is the Casimir element and
$\delta:U(\Gg)\to U(\Gg)\otimes U(\Gg)$ is
the comultiplication. Furthermore
$C=\sum_{i=1}^{n}e_{ii}^{2}+\sum_{1\leq i<j\leq n}(e_{ii}-e_{jj})
+\sum_{1\leq i<j\leq n}e_{ji}e_{ij}$ acts on the
Verma module $M(\lambda)$ as multiplication by
$
b_{\lambda}=\sum_{i=1}^{n}\lambda_{i}^{2}+\sum_{1\leq i<j\leq n}
(\lambda_{i}-\lambda_{j}).
$
It follows that $\Omega$ stabilizes any $\Gg$-submodule of $V\otimes
M(\lambda)=L(\varepsilon_{1})\otimes M(\lambda)$ 
and that the induced action on
any subquotient isomorphic to $M(\lambda+\varepsilon_{i})$ is as
multiplication by $\frac{1}{2}(b_{\lambda+\varepsilon_{i}}
-b_{\varepsilon_{1}}-b_{\lambda})=\lambda_{i}-i+1$.
Since $V\otimes M(\lambda)=\bigoplus_{a\in\BC}E_{a}M(\lambda)$,
this identifies $E_{a}M(\lambda)$ as the generalised
$a$-eigenspace of $X_{M(\lambda)}$. We deduce that for any $M\in\CO$,
the generalized $a$-eigenspace of $X_{M}$ is $E_{a}M$.

\begin{rem}
The canonical adjunction between $V\otimes-$ and $V^{\ast}\otimes-$ is
given by the canonical maps
$\eta:\BC\to V^{\ast}\otimes V$ and
$\varepsilon:V\otimes V^{\ast}\to\BC, v\otimes\xi\mapsto \xi(v)$.
Let $X_{M} \in\End_{\Gg}(V^{\ast}\otimes M)$ and
$T_{M}\in\End_{\Gg}(V^{\ast}\otimes V^{\ast}\otimes M)$ be the induced
endomorphisms (cf \S \ref{dualend}).
Then $X_{M}(\varphi\otimes m)=(-\Omega-n)(\varphi\otimes m)$
and $T_{M}(\varphi\otimes\varphi^{\prime}\otimes m)=\varphi^{\prime}
\otimes\varphi\otimes m$.
\end{rem}

\subsection{Rational representations}
\label{ratrep}
\subsubsection{} The construction of $\Gsl_2$-categorifications
in \S\ref{CatO} works, more or less in the same way, on
the category $G\mMod$ of
finite-dimensional rational representations of $G=GL_n(k)$, where
 $k$ is an algebraically closed field of characteristic $p>0$.

Denote by $\CX$ the character group of the subgroup of diagonal matrices
in $G$. 
 We identify $\CX$ with $\BZ^n$ via 
the isomorphism sending $(\lambda_1,\ldots,\lambda_n)\in\BZ^n$
 to $\lambda=\sum_i\lambda_i \varepsilon_i\in\CX$, where
 $\varepsilon_i$ is defined by $\varepsilon_i(\operatorname{diag}\nolimits
  (t_1,\ldots,t_n))=t_i$.
 This identifies the set $\CX_+$ of dominant weights with
 $\{(\lambda_1,\ldots,\lambda_n)\in\BZ^n\mid \lambda_1\geq\ldots\geq\lambda_n\}$.
For each $\lambda\in\CX_+$, let $L(\lambda)$ be the unique
simple $G$-module with highest weight
$\lambda$.

Let $B$ be the Borel subgroup of upper triangular matrices in $G$. 
For each $\lambda\in\CX$, the cohomology groups $H^i(\lambda)$ of the
associated line bundle on $G/B$ are objects of $G\mMod$. The alternating sums
$\chi(\lambda)=\sum_{i\geq 0}\ch(H^i(\lambda))\in\BZ[\CX]$ span the image
of the embedding $\ch:K_0(G\mMod)\to\BZ[\CX]$.

The Weyl group $W=\GS_n$ of $G$ acts on $\CX=\BZ^n$ by place permutations. This
extends to an action of the affine Weyl group $W_p$ generated by $W$ together
with the translations by $p\varepsilon_i-p\varepsilon_{i+1}, 1\leq i \leq n-1$.
Let $Y$ be the group of permutations of $\BZ$ generated by
$d, \sigma_0, \ldots, \sigma_{p-1}$, where $md=m+1$
and
$$m\sigma_a=\begin{cases}
m+1 & \text{ if } m\equiv a-1 \pmod p \\
m-1 & \text{ if } m\equiv a \pmod p \\
m & \text { otherwise.}
\end{cases}$$
The action of $W_p$ on $\CX=\BZ^n$ commutes with the diagonal action
of $Y$.

\begin{lemma}\label{Yorbits}
Two elements $\lambda, \mu\in\CX$ have the same stabilizer in $W_p$ if and only
of they are in the same $Y$-orbit.
\end{lemma}

\begin{proof}
Both conditions are equivalent to the following: for all $i$, $j$, and
$r$, we have $\lambda_i-\lambda_j=pr$ if and only if $\mu_i-\mu_j=pr$.
\end{proof}

We shall use the corresponding `dot actions' obtained by conjugating by
the translation by $\rho=(0,-1,\ldots,-n+1)\in\CX$:
$$w\cdot\lambda=w(\lambda+\rho)-\rho,\qquad \lambda\cdot y=
(\lambda+\rho)y-\rho.$$

Let $\Theta$ be the set of orbits of the dot action of $W_p$ on $\CX$.
For each $\theta\in\Theta$, let $\CM_\theta$ be the full subcategory
of $G\mMod$ consisting of modules whose composition factors are all
of the form $L(\lambda)$ for $\lambda\in\theta$. The Linkage Principle
\cite{CaLu}
implies that $G\mMod$ decomposes as a direct sum $G\mMod=\bigoplus
_{\theta\in\Theta}\CM_\theta$.
Let $\pr_\theta:G\mMod\to G\mMod$ denote the projection onto $\CM_\theta$.
Given $\lambda,\mu\in\CX$ and $a\in{0,\ldots,p-1}$, we write
$\lambda\to_{a}\mu$ if there exists $j$ such that $(\lambda_j-j+1)+1=
\mu_j-j+1\equiv a \pmod p$ and $\lambda_i=\mu_i$ for $i\neq j$. Note that
$\lambda\to_{a}\mu$ implies that $w\cdot\lambda\to_{a}
w\cdot\mu$ for all $w\in W_p$.  For
$\theta,\theta'\in\Theta$, we write $\theta\to_{a}\theta'$ if there
exist $\lambda\in\theta$ and $\mu\in\theta'$ such that
$\lambda\to_{a}\mu$.

Let $V$ be the natural $n$-dimensional representation of $G$. The left and
right adjoint
functors $V\otimes -:G\mMod\to G\mMod$ and $V^*\otimes -:G\mMod\to G\mMod$
 decompose as direct sums
$\bigoplus_{0\leq a\leq p-1}E_{a}$ and $\bigoplus_{0\leq a\leq p-1}F_{a}$,
 where $E_a$ and $F_a$ are sums of translation functors, defined in 
the same way as in \S\ref{CatO}. The functors $E_a$ and $F_a$ have been studied
extensively by Brundan and Kleshchev \cite{BrKl}.

Let $e_a$ and $f_a$ be the maps on characters induced by $E_a$ and $F_a$.
For each $\lambda\in\CX$, we have (eg using \cite[Proposition 7.8]{Jan})

\begin{equation} \label{EFonChi}
e_{a}\chi(\lambda)=\sum_{\substack{\mu\in\CX\\
\lambda\to_{a}\mu}}\chi(\mu),\qquad f_{a}\chi(\lambda
)=\sum_{\substack{\mu\in\CX\\
\mu\to_{a}\lambda}}\chi(\mu)
\end{equation}
in $\BZ[\CX]$. Hence
$$
e_{a}f_{a}\chi(\lambda)-f_{a}e_{a}\chi(\lambda)=c_{\lambda,a}
\chi(\lambda),
$$
where
$c_{\lambda,a}=\#\{i\mid\lambda_{i}-i+1\equiv a \pmod p\}
-\#\{i\mid\lambda_{i}-i+1\equiv a-1 \pmod p\}$.
We deduce that for each $a\in\{0,\ldots,p-1\}$
 the functors $E_{a}$ and $F_{a}$
give a weak $\Gsl_2$-categorification in which the simple module
 $L(\lambda)$ has
weight $c_{\lambda,a}$. 

\subsubsection{}
These weak $\Gsl_2$-categorifications can be improved 
to $\Gsl_2$-categorifications using the same
procedure as in the characteristic zero case \S\ref{CatO}.
We first define endomorphisms $X$ of
$V\otimes -$ and $T$ of $V\otimes V\otimes -$. Note that
 to define 
$X$, we first pass from $G$-modules to modules over
$\operatorname{Lie}\nolimits
(G)=\Ggl_n(k)$. One small modification to the argument is
required when $p$=2: in order to identify $E_a$ with the
generalized $a$-eigenspace of $X$, we write $\Omega=
-\delta(Z_2)+1\otimes Z_2 + Z_2\otimes 1 + Z_1 \otimes Z_1
 - \frac{n(n+1)}{2}$,
where $Z_1=\sum_{1\leq i\leq n}e_{ii}$ and
$Z_2=\sum_{1\leq i<j \leq n}(e_{ii}-i)(e_{jj}-j)
-\sum_{1\leq i<j \leq n} e_{ji}e_{ij}$
are central elements of
 $\operatorname{Dist}\nolimits(G)$ (cf \cite[\S 2.2]{CaLu}).

By composing the
derived (and homotopy)
 equivalences arising from these $\Gsl_2$-categorifications on
$G\mMod$, we obtain many
equivalences.

\begin{thm}\label{JeremyConj} 
Let $\lambda$ and $\mu$ be any two weights in $\CX$ with the same stabilizer
under the dot action of $W_p$. Then there are equivalences
$$K^b(\CM_{W_p\cdot\lambda})\iso K^b(\CM_{W_p\cdot\mu})\ \text{ and }\
D^b(\CM_{W_p\cdot\lambda})\iso D^b(\CM_{W_p\cdot\mu})$$
that induce the map
$$\chi(w\cdot\lambda)\mapsto\chi(w\cdot\mu)$$
on characters.
\end{thm}

\begin{rem}
Rickard conjectured the existence of such
equivalences for any
connected reductive group having a
simply-connected derived subgroup and whose root system has 
Coxeter number $h<p$ \cite[Conjecture 4.1]{Ri2}. He proved the truth
 of his conjecture in the
case of trivial stabilizers (under the weaker assumption $h\leq p$).
We do not place any restriction on $p$ in Theorem~
\ref{JeremyConj}.
\end{rem}

\begin{proof}
 By Lemma~\ref{Yorbits} we may assume that $\mu=\lambda\cdot y$ where
$y\in\{d,\sigma_0,\ldots,\sigma_{p-1}\}$. If $\mu=\lambda\cdot d$, then 
we have an equivalence $L(1,\ldots,1)\otimes-:
\CM_{W_p\cdot\lambda}\iso\CM_{W_p\cdot\mu}$,
given by tensoring with the determinant representation, that induces 
the desired map on characters.

Suppose that $\mu=\lambda\cdot\sigma_a$. Using the
$\Gsl_2$-categorification on $G\mMod$ provided by $E=E_a$ and $F=F_a$, we obtain
a self-equivalence $\Theta$ of 
$K^b(G\mMod)$ and of $D^b(G\mMod)$
such that $[\Theta]=s$ (Theorem~\ref{mainthm}).
We define an $\Gsl_2$-module $U=\bigoplus_{i\in\BZ}\BZ u_i$ by 
$eu_i=u_{i+1}$ for $i\equiv a-1 \pmod p$ and $eu_i=0$ otherwise,
and $fu_i=u_{i-1}$ for $i\equiv a \pmod p$ and $fu_i=0$ otherwise.
Then $su_i=u_{i+1}$ if $i\equiv a-1 \pmod p$, 
$su_i=-u_{i-1}$ if $i\equiv a \pmod p$,
and $su_i=u_i$ otherwise.
Thus on the tensor power $U^{\otimes n}$ we have
$su_\nu=(-1)^{h_-(\nu)}u_{\nu\sigma_a}$, where
$u_\nu=u_{\nu_1}\otimes\cdots\otimes u_{\nu_n}$
and $h_-(\nu)=\#\{i\mid\nu_{i}\equiv a \pmod p\}$.

By (\ref{EFonChi}) we have a homomorphism of
$\Gsl_2$-modules $U^{\otimes n}\to K_0
(G\mMod)$,
$u_{\nu+\rho}\mapsto\chi(\nu)$.
It follows that $s\chi(\nu)=(-1)^{h_-(\nu+\rho)}\chi(\nu\cdot\sigma_a)$.
Hence $s\chi(w\cdot\lambda)=(-1)^{h_-}\chi(w\cdot\mu)$,
where $h_-=h_-(w\cdot\lambda+\rho)=h_-(\lambda+\rho)$.
We conclude that 
$\Theta[-h_-]$ restricts to equivalences
$K^b(\CM_{W_p\cdot\lambda})\iso K^b(\CM_{W_p\cdot\mu})$
and
$D^b(\CM_{W_p\cdot\lambda})\iso D^b(\CM_{W_p\cdot\mu})$
that induce the desired map
on characters.
\end{proof}

\subsection{$q$-Schur algebras}
\label{qSchur}
We explain in this part how to obtain $\Gsl_2$-categorifications, and hence
derived equivalences, for $q$-Schur algebras.

Let $q\in k^\times$. Let
$Y_n=\bigoplus_\lambda \Ind_{H_\lambda^f}^{H_n^f}k$, where
$\lambda=(\lambda_1\ge\cdots\ge\lambda_r)$ runs over the partitions of $n$ and
$H_\lambda^f=\bigoplus_{w\in \GS_{[1,\lambda_1]}\times\cdots\times
\GS_{[n-\lambda_r+1,n]}}T_w k$ is the corresponding parabolic subalgebra
of $H_n^f$ and $k$ corresponds to the representation $1$.
We define the $q$-Schur algebra $S_n=\End_{H_n^f}(Y_n)$.

Let $\CY_n$ be the full subcategory of $H_n^f\mMod$ whose objects are
direct sums of direct summands of $Y_n$ (``$q$-Young modules'') and
let $\CY=\bigoplus_{n\ge 0}\CY_n$. Mackey's formula shows that
$\CY$ is stable under $E$ and $F$.
For each of the $\Gsl_2$-categorifications on 
$\bigoplus_{n\ge 0}H_n^f\mMod$ constructed in \S\ref{cyclo}
we deduce from Theorem \ref{restrictioncat} an $\Gsl_2$-categorification
on $\bigoplus_{n\ge 0}S_n\mMod$ and a morphism of $\Gsl_2$-categorifications
$\bigoplus_{n\ge 0}S_n\mMod\to \bigoplus_{n\ge 0}H_n^f\mMod$. This
provides a version of Theorem \ref{equivHecke} for $q$-Schur algebras.

\begin{rem}
We go back to the setting of \S\ref{GLq} (in particular,
$q$ is a prime power). The canonical map
$A_n\to \End_{H_n^{f\opp}}(kG_n e_{B_n})^\opp$ is surjective and
its image $S'_n$ is Morita equivalent to $S_n$
(``double centralizer Theorem'', cf \cite{Ta}). This gives by restriction
a fully faithful functor $S_n\mMod\iso S'_n\mMod\to A_n\mMod$. Since
$E(kG_n e_{B_n})\simeq kG_{n+1}e_{B_{n+1}}$, it follows that
$\bigoplus_{n\ge 0}S_n\mMod$ is stable under $E$. Mackey's formula
shows that it is also stable under $F$. This gives a morphism
of weak $\Gsl_2$-categorifications $\bigoplus_{n\ge 0}S_n\mMod\to
\bigoplus_{n\ge 0}A_n\mMod$ and the composition with the
morphism $\bigoplus_{n\ge 0}A_n\mMod\to
\bigoplus_{n\ge 0}H_n^f\mMod$ of \S\ref{GLqcat} is isomorphic to the morphism
$\bigoplus_{n\ge 0}S_n\mMod\to \bigoplus_{n\ge 0}H_n^f\mMod$ constructed above.
One deduces that  $\bigoplus_{n\ge 0}S_n\mMod\to
\bigoplus_{n\ge 0}A_n\mMod$ is actually a morphism of
$\Gsl_2$-categorifications.

Note also that we get another proof of Lemma \ref{sl2relGL} using the fact
that
the canonical map $K_0(S_n\mMod)\iso K_0(A_n\mMod)$ is an isomorphism.
\end{rem}

\begin{rem}
The interested reader will extend the results of \S\ref{ratrep} to the quantum
case and show that the categorification of $q$-Schur algebras can
be realized as a subcategorification of the quantum group case.
\end{rem}

\subsection{Realizations of minimal categorifications}
\label{realmin}
\subsubsection{}
We now show that the minimal categorification of \S\ref{mincat} is a special
case of the categorification on representations of blocks 
of cyclotomic Hecke algebras.

\smallskip
Fix $a\in k^\times$ and put $v=(v_1,\ldots,v_n)=(a,\ldots,a)$. 
Then  
$\CH_i=\CH_i(q,v)$ is the quotient of $H_i$ be the ideal
generated by $x_1^n$ (where $x_1=X_1-a$).
The kernel of the action of $H_i$ on the simple module
$K_i=H_i\otimes_{P_i}P_i/\Gm_i$
contains $x_1^n$ if and only if $i\leq n$ (cf \S\ref{en});
let $A_i$ be the block of $\CH_i$ containing $K_i$ for $0\leq i \leq n$.
A finitely generated $\CH_i$-module $M$ is in $A_i$ if and only if
$\Gn_i$ acts nilpotently on $M$ (equivalently $\Gm_i$ acts nilpotently 
on $M$), and $K_i$ is the unique simple module in $A_i$.

We have $FM=0$ for $M\in A_0\mMod$ and
$FM=\Res_{\CH_{i-1}}^{\CH_i}M\in A_{i-1}\mMod$ for
$M\in A_i\mMod$ and $0<i\leq n$.

Likewise $EM=0$ for $M\in A_n\mMod$.
Let $M\in A_i\mMod$ with $0\leq i<n$.
Consider $N$ a simple $\CH_{i+1}$-quotient of $EM$.
We have $\Hom(EM,N)\simeq \Hom(M,FN)\not=0$. In particular,
$FN$ has a non-zero $\CH_i$-submodule $M'$ on which $x_1,\ldots,x_i$
act nilpotently. Let $M''$ be the $(k[x_{i+1}]\otimes \CH_i)$-submodule
of $FN$ generated by $M'$. Then, $x_1,\ldots,x_{i+1}$ act nilpotently
on $M''$. Now, $N$ is a simple $\CH_{i+1}$-module, hence it is generated
by $M''$ as a $\CH_{i+1}$-module, so  $x_1,\ldots,x_{i+1}$ act nilpotently
on $N$. We deduce that they act nilpotently on $EM$ as well.
So, $EM\in A_{i+1}\mMod$.

\smallskip
So $\CA=\bigoplus_i A_i\mMod$ is an ${\Gsl}_2$-categorification and
$\BQ\otimes K_0(\CA)$ is a simple $\Gsl_2$-module of dimension $n+1$.
Let $U=K_0=k$, the simple (projective) module for $A_0=k$. The
morphism of ${\Gsl}_2$-categorifications
$R_U:\CA(n)\to\CA$ is an equivalence (Proposition \ref{mincatequiv}).
In particular $\bar{H}_{i,n}$ and $A_i$ are isomorphic, as each has
an $i!$-dimensional simple module.

\subsubsection{}
We explained in \S\ref{grassmann} that $\bar{H}_{i,n}$ is Morita equivalent
to its center, which is isomorphic to the cohomology of certain
Grasmmannian varieties. We sketch here a realization
of the minimal categorification in that
setting. We consider only the case $q=1$; the case $q\not=1$ can be
dealt with similarly, replacing cohomology by $\BG_m$-equivariant
$K$-theory.

Let $G_{i,j}$ be the variety of pairs $(V_1,V_2)$ of subspaces
of $\BC^n$ with $V_1\subset V_2$, $\dim V_1=i$ and $\dim V_2=j$.
We put $A_i=H^*(G_i)$. The $(A_{i+1},A_i)$-bimodule $H^*(G_{i,i+1})$
defines by tensor product
a functor $E_{i}:A_i\mMod\to A_{i+1}\mMod$ and switching sides,
a left and right adjoint $F_{i}:A_{i+1}\mMod\to A_i\mMod$.
Let $E=\bigoplus E_{i}$ and $F=\bigoplus F_{i}$.
This gives a weak $\Gsl_2$-categorification that has been considered
by Khovanov as a way of categorifying irreducible $\Gsl_2$-representations.
It is a special case of the construction of irreducible
finite dimensional representations of $\Gsl_n$ due to Ginzburg \cite{Gi}.

\smallskip
We denote by $X$ the endomorphism of $E$ given on $H^*(G_{i,i+1})$
by cup product by $c_1(L_{i+1})$.
We have a $\BP^1$-fibration $\pi:G_{i,i+1}\times_{G_{i+1}}G_{i+1,i+2}\to
G_{i,i+2}$ given by first and last projection. It induces a
structure of $H^*(G_{i,i+2})$-module on
$H^*(G_{i,i+1}\times_{G_{i+1}}G_{i+1,i+2})
=H^*(G_{i,i+1})\otimes_{H^*(G_{i+1})}H^*(G_{i+1,i+2})$. There is
a unique endomorphism $T$ of $H^*(G_{i,i+2})$-module on
$H^*(G_{i,i+1}\times_{G_{i+1}}G_{i+1,i+2})$
satisfying $T(c_1(L_{i+1}))=c_1(L_{i+2})-1$. This provides us with an
endomorphism of $E_{i+1}E_{i}$ and taking the sum over all $i$, we
get an endomorphism $T$ of $E^2$.
One checks easily that this gives an $\Gsl_2$-categorification
(with $a=0$) that is isomorphic to the minimal categorification.

\smallskip
The functor $E^{(1,r)}:A_i\mMod\to A_{i+r}\mMod$ is isomorphic to
the functor given by the bimodule $H^*(G_{i,i+r})$.

Take $i\le n/2$ and let us
now consider $\Theta[-i]$, restricted to a functor
$D^b(H^*(G_i)\mMod)\iso D^b(H^*(G_{n-i})\mMod)$. It is probably isomorphic
to the functor given by the cohomology of the subvariety
$\{(V,V')|V\cap V'=0\}$ of $G_i\times G_{n-i}$, the usual kernel
for the Grassmannian duality (cf eg \cite[Exercice III.15]{KaScha}).

\end{document}